\documentclass[10pt]{article}

\RequirePackage{amsmath}
\RequirePackage{amsfonts}
\RequirePackage{amssymb}
\RequirePackage{hyperref}
\RequirePackage{algorithm}
\RequirePackage{algpseudocode}
\RequirePackage{graphicx}
\RequirePackage{xcolor}
\RequirePackage{amsthm}
\allowdisplaybreaks
\hypersetup{colorlinks, hypertexnames=false, pageanchor=true,
            linkcolor=blue, citecolor=violet, urlcolor=magenta,
            pdftitle={stoc-ipm}
            }

\usepackage{dpr_fec,ifthen,dsfont}
\usepackage{tikz}

\newtheorem{parameter}{Parameter}
\newtheorem{example}{Example}

\usepackage{isetechreport}
\def\reportyear{24T}
\def\reportno{008}
\def\revisionno{0}
\def\originaldate{\today}
\coraltrue
\cvcrfalse

\begin{document}
\title{
Single-Loop Deterministic and Stochastic Interior-Point Algorithms for Nonlinearly Constrained Optimization

}

\author[1]{Frank E.~Curtis\thanks{E-mail: frank.e.curtis@lehigh.edu}}
\author[1]{Xin Jiang\thanks{E-mail: xjiang@lehigh.edu}}
\author[1]{Qi Wang\thanks{E-mail: qiw420@lehigh.edu}}
\affil[1]{Department of Industrial and Systems Engineering, Lehigh University}
\titlepage

\maketitle

%**********
% Abstract
%**********
\begin{abstract}
An interior-point algorithm framework is proposed, analyzed, and tested for solving nonlinearly constrained continuous optimization problems.  The main setting of interest is when the objective and constraint functions may be nonlinear and/or nonconvex, and when constraint values and derivatives are tractable to compute, but objective function values and derivatives can only be estimated.  The algorithm is intended primarily for a setting that is similar for stochastic-gradient methods for unconstrained optimization, namely, the setting when stochastic-gradient estimates are available and employed in place of gradients of the objective, and when no objective function values (nor estimates of them) are employed.  This is achieved by the interior-point framework having a single-loop structure rather than the nested-loop structure that is typical of contemporary interior-point methods.  For completeness, convergence guarantees for the framework are provided both for deterministic and stochastic settings.  Numerical experiments show that the algorithm yields good performance on a large set of test problems.
\end{abstract}

%**************
% Body of paper
%**************
%************************
% Paper-dependent macros
%************************
\newcommand{\barrier}{\tilde\phi}
\newcommand{\buff}{\mathrm{buff}}

%*********
% Section
%*********
\section{Introduction}\label{sec.introduction}

The interior-point methodology is one of the most popular and effective derivative-based algorithm methodologies for solving inequality-constrained continuous optimization problems.  Indeed, while there has been a long history of development of other types of algorithms---e.g., of the augmented-Lagrangian and sequential-quadratic-optimization varieties \cite{NW99}---interior-point methods are arguably the most effective class of approaches for solving many constrained continuous optimization problems.  In the setting of nonconvex optimization, which is the setting of interest in this paper, the effectiveness of interior-point methods is evidenced by the fact that the most popular software packages for solving such problems are based on interior-point methods.  These include the packages \texttt{Ipopt} \cite{WaecBieg06}, \texttt{Knitro} \cite{ByrdHribNoce99}, and \texttt{LOQO} \cite{VandShan99}.  This popularity withstands despite the fact that interior-point methods for solving nonconvex problems do not necessarily achieve the same optimal worst-case complexity properties that can be achieved in convex settings.

Despite its widespread success, little work has been done on extending the interior-point methodology to the stochastic regime.  Here, we are focused on the setting in which constraint values and derivatives are tractable to compute, but objective values and derivatives are not.  Furthermore, like in the setting of stochastic-gradient-based methods for solving unconstrained problems, our focus is to have an interior-point method that employs stochastic-gradient estimates (in place of exact objective gradients) and no objective function evaluations.  Given the success of interior-point methods in the deterministic regime, there is ample motivation to explore their potential theoretical and practical behavior in such a stochastic regime.  This is the main motivation for the work presented in this paper.

We contend that there are at least a few significant obstacles to the design and analysis of a stochastic interior-point method, which may be at least part of the reason that there has been little prior work on this topic.  Firstly, interior-point methods often have a nested-loop structure wherein an inner loop is designed to solve (approximately) a so-called barrier subproblem for a fixed value of a barrier parameter and the outer loop involves decreasing the barrier parameter.  Determining when to terminate the inner loop typically involves evaluating a measure of stationarity.  However, in the aforementioned stochastic setting, evaluating such a measure is not tractable since objective-gradient evaluations are not tractable.  Secondly, the barrier functions introduced in interior-point methods have corresponding gradient functions that are not Lipschitz continuous nor bounded.  This presents major difficulties since for many stochastic-gradient-based methods, Lipschitz continuity and boundedness of the gradients are necessary features for ensuring convergence guarantees.  Thirdly, interior-point methods regularly involve globalization mechanisms such as line searches or trust-region step-acceptance mechanisms that require exact function evaluations, but these are not tractable to perform in the aforementioned stochastic setting that is of interest here.

%************
% Subsection
%************
\subsection{Contributions and Limitations}\label{sec.contributions}

To provide a significant advancement in the design and analysis of stochastic-gradient-based interior-point methods, in this paper we propose and analyze a single-loop interior-point algorithm framework that offers rigorous convergence guarantees in both deterministic and stochastic settings.  Let us emphasize upfront that we do not expect our method to be competitive with state-of-the-art interior-point methods in the deterministic regime, and especially not with those that employ second-order derivatives of the objective.  Nevertheless, we present an algorithm and analyze its convergence properties for that setting for completeness and since it helps to elucidate the essential features of an algorithm that are necessary for achieving convergence guarantees in the stochastic regime.

Our method is an extension of the approach proposed in \cite{CKRW23}, which was proposed and analyzed for the bound-constrained setting.  We go beyond the bound-constrained setting and propose a framework that is applicable for problems with nonlinear inequality constraints, at least as long as a strictly feasible initial point for the algorithm can be provided.  For the sake of generality, we present the framework in a manner such that linear equality constraints can also be present, although we emphasize that our convergence guarantees in neither the deterministic nor stochastic regime allow for the presence of nonlinear equality constraints.  Along with our concluding remarks, we discuss the challenges that would be faced when trying to extend our approach to settings with nonlinear equality constraints.

Most of the assumptions that are required for our convergence guarantees are standard in the literature on interior-point algorithms, specifically for that on so-called feasible interior-point methods that maintain feasibility (with respect to all of the constraints) at all algorithm iterates.  However, we require one additional assumption that is not standard.  This assumption essentially requires that, when an algorithm iterate is close to a constraint boundary, a direction of sufficient descent can be computed that moves sufficiently away from the boundary.  It is, in essence, an assumption that combines two aspects: (a) a type of nondegeneracy assumption and (b) an assumption that the barrier parameter is initialized to be sufficiently large compared to a neighborhood parameter that is introduced in our algorithm.  We demonstrate these aspects through some illustrative examples.  Since, to the best of our knowledge, there are no other interior-point methods in the literature that possess convergence guarantees in a stochastic setting that is comparable to ours, we contend that our algorithm and analysis together constitute a notable contribution to the literature despite our method being a feasible interior-point method and our need to introduce a nonstandard assumption.  Future research may reveal techniques that can offer convergence guarantees for a stochastic-gradient-based interior-point method in other settings.

We also provide the results of numerical experiments showing that our algorithms performs well when employed to solve challenging test problems.

%************
% Subsection
%************
\subsection{Notation}\label{sec.notation}

We use $\R{}$ to denote the set of real numbers and, given $a \in \R{}$, we use $\R{}_{\geq a}$ (resp.,~$\R{}_{>a}$) to denote the set of real numbers that are greater than or equal to $a$ (resp.,~strictly greater than $a$).  Superscripts are used to indicate the dimension of a vector or matrix whose elements are restricted to such a set; e.g., the set of $n$-dimensional vectors is denoted as $\R{n}$.  The set of positive integers is denoted as $\N{} := \{1,2,\dots\}$ and, given $n \in \N{}$, the set of positive integers up to $n$ is denoted as $[n] := \{1,\dots,n\}$.  The set of $n$-by-$n$-dimensional symmetric positive semidefinite (resp.,~definite) matrices is denoted as $\Smbb_{\succeq0}^n \subset \R{n \times n}$ (resp.,~$\Smbb_{\succ0}^n \subset \R{n \times n}$).

A real-valued sequence is always introduced as $\{v_k\} \subset \Rcal$, where $\Rcal$ is a real vector (sub)space and the subset notation indicates that $v_k \in \Rcal$ for each index $k \in \N{}$ (although elements in the sequence may repeat).  Here, subscripts are used to indicate index number in a sequence, although in some situations a subscript is used to indicate an element number of a vector or vector-valued function.  When indicating the $i$th element of a vector $v_k$, the notation $[v_k]_i$ is used.  In all cases, the meaning of a subscript is clear from the context.  The notation $\{v_k\} \to v$ indicates that the sequence $\{v_k\}$ converges to $v$.  The notation $\{v_k\} \searrow 0$ indicates that $\{v_k\}$ is a sequence of positive real numbers that vanishes monotonically.

Given a vector $v \in \R{n}$, we use $\diag(v)$ to denote the diagonal matrix for which, for all $i \in [n]$, the $(i,i)$ element is $v_i \in \R{}$. Given two vectors $z \in \R{m}$ and $c \in \R{m}$, their Hadamard (i.e., component-wise) product is denoted as $z \circ c \in \R{m}$.  We use $\ones$ to denote a vector whose elements are all equal to one; the length of such a vector is determined by the context in which it appears.  Given $A \in \R{l \times n}$, the null space of $A$ is denoted as $\Null(A)$.  Given $(A,B) \in \Smbb_{\succeq0}^n \times \Smbb_{\succeq0}^n$, we use $A \succeq B$ (resp., $A \preceq B$) to indicate that $A - B$ (resp., $B - A$) is positive semidefinite.  Given $v \in \R{p}$ and $M \in \Smbb_{\succ 0}^p$, the $M$-norm of $v$ is written as $\|v\|_M = \sqrt{v^TMv}$.  Finally, given a set $\Pcal$, its interior is denoted as $\interior(\Pcal)$ and its polar cone is denoted as
\bequationNN
  \Pcal^\circ := \{ y \in \R{n} : p^Ty \leq 0\ \text{for all}\ p \in \Pcal \}.
\eequationNN

%************
% Subsection
%************
\subsection{Outline}\label{sec.outline}

In \S\ref{sec.problem}, we present our problem of interest and preliminary commentary.  In \S\ref{sec.algorithm}, we present our proposed framework and prove basic facts about its behavior.  Section~\ref{sec.analysis} contains our assumptions and convergence analyses for the deterministic and stochastic settings.  As previously mentioned, our analyses rely on an assumption that is not standard; we discuss this assumption in detail in~\S\ref{sec.assumption}, highlighting certain settings of interest in which it indeed holds.  Details about a practical implementation of our algorithm and resulting numerical experiments are presented in~\S\ref{sec.numerical}.  Concluding remarks are provided in \S\ref{sec.conclusion}.

%*********
% Section
%*********
\section{Problem Formulation and Preliminaries}\label{sec.problem}

Our proposed algorithm framework (Algorithm~\ref{alg.slip} on page~\pageref{alg.slip}) is designed to solve nonlinearly constrained optimization problems of the form
\bequation\label{prob.opt}
  \min_{x \in \R{n}}\ f(x)\ \subjectto\ Ax = b\ \text{and}\ c(x) \leq 0,
\eequation
where $A \in \R{l \times n}$ (with $l < n$), $b \in \R{l}$, and the functions $f$ and $c$ satisfy the following assumption with respect to the feasible region of problem~\eqref{prob.opt}, namely,
\bequationNN
  \Fcal := \Ecal \cap \Ccal_{\leq0},\ \text{where}\ \left\{ \baligned \Ecal &:= \{x \in \R{n} : Ax = b\}; \\ \Ccal_{\leq0} &:= \{x \in \R{n} : c(x) \leq 0\}. \ealigned \right.
\eequationNN

\bassumption\label{ass.main}
  There exists an open set $\overline\Fcal \supset \Fcal$ over which the objective function $f : \overline\Fcal \to \R{}$ is continuously differentiable and bounded below by $f_{\inf} \in \R{}$ and the constraint function $c : \overline\Fcal \to \R{m}$ is continuously differentiable.  In addition, the matrix $A \in \R{l \times n}$ has full row rank.  Finally, there exists $x \in \Fcal$ such that $c(x) < 0$ $($i.e., where the inequality holds strictly, component-wise$)$.
\eassumption

The elements of Assumption~\ref{ass.main} pertaining to $f$ and $c$ are mostly common for the literature on smooth nonlinearly constrained optimization.  The only exception is the assumption that there exists $x \in \Fcal$ such that $c(x) < 0$.  This aspect of the assumption precludes the presence of nonlinear equality constraints (through two-sided inequalities).  It is needed in our setting since our algorithm requires a feasible initial point that is strictly feasible with respect to the nonaffine constraints.

Assumption~\ref{ass.main} does not require that the feasible region~$\Fcal$ is bounded, although our convergence analysis in~\S\ref{sec.analysis} requires that the iterates remain within a set such that $c$ remains component-wise bounded below, and that the functions $f$, $c$, $\nabla f : \overline\Fcal \to \R{n}$, and $\nabla c_i : \overline\Fcal \to \R{n}$ for all $i \in [m]$ are Lipschitz continuous; see Assumption~\ref{ass.bound} on page~\pageref{ass.bound}.  Our analysis could be extended to the setting in which the matrix~$A$ does not have full row rank as long as the linear system $Ax=b$ is consistent.  Our algorithm framework requires that the initial iterate is contained in $\Ecal$ and that each step lies in the null space of $A$, which inductively implies that each iterate lies in $\Ecal$.  These requirements remain reasonable when~$A$ has linearly dependent rows as long as a projection operator onto the null space of $A$ is available. Merely for ease of notation, we assume $A$ has full row rank.

In fact, our algorithm framework ensures that all iterates remain in the set
\bequation\label{eq.F_strict}
  \Fcal_{<0} := \Ecal \cap \Ccal_{<0},\ \text{where}\ \Ccal_{<0} := \{x \in \R{n} : c(x) < 0\}.
\eequation
Observe that this set is guaranteed to be nonempty under Assumption~\ref{ass.main}.  In each iteration, our framework involves the computation of a search direction that, for a given barrier parameter $\mu \in \R{}_{>0}$, aims to reduce the value of the objective augmented with the logarithmic barrier function defined with respect to the inequality-constraint functions of \eqref{prob.opt}, namely, the function $\phi(\cdot,\mu) : \Ccal_{<0} \to \R{}$ defined by
\bequationNN
  \phi(x,\mu) = f(x) - \mu \sum_{i\in[m]} \log (-c_i(x)).
\eequationNN
Going forward, we refer to this as the \emph{barrier-augmented objective function}.  Observe that if $x^\mu$ is a local solution of the problem to minimize $\phi(\cdot,\mu)$ over~$\Fcal_{<0} = \Ecal \cap \Ccal_{<0}$, then there exists a Lagrange multiplier vector $y^\mu \in \R{l}$ such that, along with $z_i^\mu := -\mu (c_i(x^\mu))^{-1}$ for all $i \in [m]$, one finds that
\bequationNN
  \nabla f(x^\mu) + A^Ty^\mu + \nabla c(x^\mu) z^\mu = 0,\ Ax^\mu = b,\ c(x^\mu) < 0,\ z^\mu > 0,\ \text{and}\ -z^\mu \circ c(x^\mu) = \mu \ones.
\eequationNN
Any tuple $(x^\mu,y^\mu,z^\mu)$ satisfying this system is referred to as a stationary point for the problem to minimize $\phi(\cdot,\mu)$ over $\Fcal_{<0}$.  This system should be compared to the Karush-Kuhn-Tucker (KKT) conditions for problem~\eqref{prob.opt}; in particular, under a constraint qualification---such as the Mangasarian-Fromovitz constraint qualification (MFCQ) \cite{MangFrom67}---if $x$ is a local solution of problem~\eqref{prob.opt}, then there exists $y \in \R{l}$ and $z \in \R{m}$ such that the following system of KKT conditions is satisfied:
\bequation\label{eq.KKT}
  \nabla f(x) + A^Ty + \nabla c(x) z = 0,\ Ax = b,\ c(x) \leq 0,\ z \geq 0,\ \text{and}\ -z \circ c(x) = 0.
\eequation
Any tuple $(x,y,z)$ satisfying this system is referred to as a stationary point (or KKT point) for problem~\eqref{prob.opt}.  As is generally the case for an interior-point method, our algorithm framework aims to find such a KKT point by taking steps to reduce the barrier-augmented objective function for a diminishing sequence of barrier parameters---a procedure that can be motivated by the fact that the stationarity conditions for the minimization of $\phi(\cdot,\mu)$ over $\Fcal_{<0}$ tend toward those of~\eqref{prob.opt}.

%*********
% Section
%*********
\section{Algorithm}\label{sec.algorithm}

Our algorithm framework is stated as Algorithm~\ref{alg.slip} on page~\pageref{alg.slip}.  The framework presumes that it is given an initial neighborhood parameter $\theta_0 \in \R{}_{>0}$ and an initial iterate $x_1 \in \Ecal \cap \Ncal(\theta_0)$, where for all $\theta \in \R{}_{>0}$ we define
\bequationNN
  \Ncal(\theta) := \{x \in \R{n} : c(x) \leq -\theta \ones\} \subseteq \Ccal_{<0}.
\eequationNN
Observe that the existence of such $\theta_0$ and $x_1$ is guaranteed under Assumption~\ref{ass.main}.  Moreover, the steps of Algorithm~\ref{alg.slip} ensure that, for a prescribed, monotonically decreasing sequence $\{\theta_k\} \subset \R{}_{>0}$ with $\theta_1 < \theta_0$, one finds $x_k \in \Ecal \cap \Ncal(\theta_{k-1})$ for all $k \in \N{}$.  Hence, to describe the main steps of the algorithm for each $k \in \N{}$, we may presume that for all $k \in \N{}$ the algorithm has an iterate $x_k \in \Ecal \cap \Ncal(\theta_{k-1})$.

Following \cite{CKRW23}, Algorithm~\ref{alg.slip} also employs a prescribed, monotonically decreasing sequence of barrier parameters $\{\mu_k\} \subset \R{}_{>0}$ such that $\{\mu_k/\theta_{k-1}\}$ is a constant sequence.  Other inputs required by Algorithm~\ref{alg.slip} in the deterministic setting are $\eta \in (\theta_0/\mu_1,1)$, $\underline\eta \in (\theta_0,\infty)$, $\underline\zeta \in \R{}_{>0}$, $\overline\zeta \in \R{}_{>0}$ with $\underline\zeta \leq \overline\zeta$, and $\zeta \in (0,1)$, each of which influences the search direction computation.  Given these inputs along with $k \in \N{}$ and $x_k \in \Ecal \cap \Ncal(\theta_{k-1})$, the computation proceeds as follows.  Firstly, a set of indices of constraints that are \emph{nearly active} with respect to $x_k$ is identified as
\bequationNN
  \Acal_k := \{i \in [m] : -\eta \mu_k < c_i(x_k)\}.
\eequationNN
Observe that since $x_k \in \Ncal(\theta_{k-1})$ it follows that $c_i(x_k) \leq -\theta_{k-1}$ for all $i \in [m]$, and, since the facts that $\eta \in (\theta_0/\mu_1,1)$ and $\{\mu_k/\theta_{k-1}\}$ is a constant sequence together ensure $-\eta \mu_k < -\theta_{k-1}$, it follows that $-\eta \mu_k < c_i(x_k) \leq -\theta_{k-1}$ for all $i \in \Acal_k$.  Secondly, a vector $q_k \in \R{n}$ is computed, the value of which is distinct between the deterministic and stochastic settings.  In particular, the framework employs
\bequation\label{eq.gradient}
  q_k := \bcases \nabla f(x_k) - \mu_k \nabla c(x_k) \diag(c(x_k))^{-1}\ones & \text{(deterministic)} \\ g_k - \mu_k \nabla c(x_k) \diag(c(x_k))^{-1}\ones & \text{(stochastic)} \ecases
\eequation
where in the stochastic setting the vector $g_k$ is a stochastic estimate of $\nabla f(x_k)$.  Third, with $P := I - A^T(AA^T)^{-1}A$ denoting the orthogonal projection matrix onto the null space of $A$, a search direction $d_k \in \R{n}$ is computed such that
\bsubequations\label{eq.d_conditions}
  \begin{align}
    Ad_k &= 0, \label{eq.d_condition_1} \\
    \underline\zeta \|Pq_k\|_2 &\leq \|d_k\|_2, \label{eq.d_condition_2} \\
    \overline\zeta \|Pq_k\|_2 &\geq \|d_k\|_2, \label{eq.d_condition_3} \\
    -(Pq_k)^Td_k &\geq \zeta \|Pq_k\|_2 \|d_k\|_2,\ \ \text{and} \label{eq.d_condition_4} \\
    \nabla c_i(x_k)^Td_k &\leq - \thalf \underline\eta \|d_k\|_2\ \ \text{for all}\ \ i \in \Acal_k. \label{eq.d_condition_5}
  \end{align}
\esubequations
In the deterministic setting, all that is required of the search direction $d_k$ is that it satisfies \eqref{eq.d_conditions}.  For the stochastic setting, we introduce in our analysis a particular strategy for computing the search direction that implies that \eqref{eq.d_conditions} holds along with other useful properties for our analysis for that setting.

Let us now discuss \eqref{eq.d_conditions} in detail.  Firstly, condition~\eqref{eq.d_condition_1} requires that the search direction satisfies $d_k \in \Null(A)$, i.e., it lies in the null space of~$A$, which, since $x_k \in \Ecal$, ensures that the subsequent iterate will also lie in $\Ecal$.  Secondly, conditions~\eqref{eq.d_condition_2} and \eqref{eq.d_condition_3} require that the norm of the search direction $d_k$ is proportional to the norm of $Pq_k = (I - A^T(AA^T)^{-1}A)q_k$, i.e., the norm of the projection of~$q_k$ onto $\Null(A)$.  Thus, conditions~\eqref{eq.d_condition_2} and~\eqref{eq.d_condition_3} require that the norm of the search direction is proportional to that of a projected gradient (estimate), which is a common type of requirement for a search direction in an algorithm for minimizing a smooth objective function over an affine set.  Thirdly, condition~\eqref{eq.d_condition_4} requires that~$d_k$ makes an angle with $Pq_k$ that is sufficiently obtuse.  Again, this is a common type of requirement; e.g., in the deterministic setting, it ensures that any nonzero $d_k$ is a direction of sufficient descent with respect to the barrier-augmented objective function. Together, conditions \eqref{eq.d_condition_1}--\eqref{eq.d_condition_4} define a nonempty set and computing a search direction to satisfy these conditions is straightforward; e.g., one can compute the projection of $-q_k$ onto $\Null(A)$ and scale the resulting direction, if needed, to satisfy \eqref{eq.d_condition_2}--\eqref{eq.d_condition_3}.  However, it remains finally to consider condition~\eqref{eq.d_condition_5}.  Unfortunately, this is a requirement that is not always satisfiable along with the other conditions, namely, \eqref{eq.d_condition_1}--\eqref{eq.d_condition_4}.  That said, we contend that all of \eqref{eq.d_conditions} can be satisfied as long as (a) the constraints are not degenerate in some sense and (b) the initial barrier parameter $\mu_1$ is sufficiently large relative to $\theta_0$.  Denoting the conical hull of nearly active constraint gradients as
\bequationNN
  \Pcal_k := \left \{ \sum_{i \in \Acal_k} \nabla c_i(x_k) \sigma_i : \sigma \in \R{|\Acal_k|}_{\geq 0} \right \},
\eequationNN
condition~\eqref{eq.d_condition_5} requires that $d_k$ lies sufficiently within the interior of the polar cone of $\Pcal_k$, i.e., $\interior(\Pcal_k^\circ)$.  Hence, to satisfy \eqref{eq.d_conditions}, there needs to exist a direction sufficiently within this interior that also satisfies \eqref{eq.d_condition_1}--\eqref{eq.d_condition_4}.  For our purposes now, let us simply introduce Assumption~\ref{ass.c} below, which implicitly means that Algorithm~\ref{alg.slip} is being applied to a problem such that a search direction satisfying \eqref{eq.d_conditions} always exists.  We defer until \S\ref{sec.assumption} a more detailed discussion Assumption~\ref{ass.c}, where in particular we provide (a) example problems for which the assumption can be shown to hold or cannot be shown to hold, (b) further explanation about why the assumption can be viewed as a type of nondegeneracy assumption and an assumption about $\mu_1/\theta_0$ being sufficient large (which is reminiscient of a requirement in Corollary~3.8 and Lemma~3.15 in \cite{CKRW23} for the bound-constrained setting), and (c) a procedure for computing $d_k$ to satisfy \eqref{eq.d_conditions} when such a direction exists.

\bassumption\label{ass.c}
  For all $k \in \N{}$ generated in a run of Algorithm~\ref{alg.slip}, there exists a search direction $d_k \in \R{n}$ satisfying the conditions \eqref{eq.d_conditions}.
\eassumption

Upon computation of the search direction, Algorithm~\ref{alg.slip} chooses a step size $\alpha_k \in \R{}_{>0}$.  For simplicity, our analyses for the deterministic and stochastic settings consider conservative rules for selecting the step sizes that each depend on the barrier- and neighborhood-parameter sequences, $\{\mu_k\}$ and $\{\theta_k\}$.  We conjecture that it would be possible to generalize our algorithm and analyses to allow more flexibility in the step-size selection rules, as is done in \cite{CKRW23}.  However, for our purposes here, we present conservative strategies that are sufficient for proving convergence.

The algorithm's last major step in the $k$th iteration is to compute a positive real number $\gamma_k \in \R{}_{>0}$ such that the line segment $[x_k,x_k + \gamma_k \alpha_k d_k]$ is contained in the neighborhood $\Ncal(\theta_k)$.  The subsequent iterate is then set as $x_{k+1} \gets x_k + \gamma_k \alpha_k d_k$.  As mentioned at the beginning of this section, from the facts that the $k$th iterate has $x_k \in \Ecal \cap \Ncal(\theta_{k-1})$, the $k$th search direction satisfies \eqref{eq.d_condition_1}, and the choice of~$\gamma_k$ ensures that $x_k + \gamma_k \alpha_k d_k \in \Ncal(\theta_k)$, it follows that $x_{k+1} \in \Ecal \cap \Ncal(\theta_k)$.  For the sake of formality, we state the following lemma that has been proved.

\blemma\label{lem.in_neighborhoods}
  For all $k \in \N{}$ generated in any run of Algorithm~\ref{alg.slip}, it follows that the iterate satisfies $x_k \in \Ecal \cap \Ncal(\theta_{k-1}) \subseteq \Fcal_{<0}$ and the search direction satifies $d_k \in \Null(A)$.
\elemma

\begin{algorithm}[ht]
  \caption{Single-Loop Interior-Point (SLIP) Method}
  \label{alg.slip}
  \begin{algorithmic}[1]
    \Require initial neighborhood parameter $\theta_0 \in \R{}_{>0}$, initial iterate $x_1 \in \Ecal \cap \Ncal(\theta_0)$, barrier-parameter sequence $\{\mu_k\} \searrow 0$, neighborhood-parameter sequence $\{\theta_k\} \searrow 0$ with $\theta_1 < \theta_0$ such that $\{\mu_k/\theta_{k-1}\}$ is a constant sequence, maximum neighborhood-parameter sequence $\{\gamma_{k,\max}\} \subset (0,1]$, $\eta \in (\theta_0/\mu_1,1)$, $\underline\eta \in (\theta_0,\infty)$, $\underline\zeta \in \R{}_{>0}$, $\overline\zeta \in \R{}_{>0}$ with $\underline\zeta \leq \overline\zeta$, and $\zeta \in (0,1)$
    \For{\textbf{all} $k\in\N{}$}
      \State Compute $q_k \in \R{n}$ by \eqref{eq.gradient}
      \State Compute $d_k \in \R{n}$ to satisfy \eqref{eq.d_conditions}
      \State Choose $\alpha_k \in \R{}_{>0}$
      \State Compute $\gamma_k \in (0, \gamma_{k,\max}]$ such that $[x_k,x_k + \gamma_k \alpha_k d_k] \subseteq \Ncal(\theta_k)$ \label{step.gamma}
      \State Set $x_{k+1} \leftarrow x_k + \gamma_k \alpha_k d_k$
    \EndFor
  \end{algorithmic}
\end{algorithm}

%*********
% Section
%*********
\section{Analysis}\label{sec.analysis}

We prove theoretical convergence guarantees for Algorithm~\ref{alg.slip} in both deterministic and stochastic settings.  For both settings, we make the following additional, standard type of assumption beyond Assumptions~\ref{ass.main} and \ref{ass.c}.

\bassumption\label{ass.bound}
  Let $\Xcal$ be an open convex set containing the sequence of line segments between iterates, namely, $\{[x_k,x_k+\alpha_kd_k]\}$, generated in any run of Algorithm~\ref{alg.slip}.  There exist constants $\kappa_{\nabla f} \in \R{}_{>0}$, $L_{\nabla f} \in \R{}_{>0}$, $\{\kappa_{c_i}\}_{i\in[m]} \subset \R{}_{>0}$, $\{L_{c_i}\}_{i\in[m]} \subset \R{}_{>0}$, $\{\kappa_{\nabla c_i}\}_{i\in[m]} \subset \R{}_{>0}$, and $\{L_{\nabla c_i}\}_{i\in[m]} \subset \R{}_{>0}$ such that:
  \bitemize
    \renewcommand\labelitemi{$\bullet$}
    \item $f$ is Lipschitz with constant $\kappa_{\nabla f}$ over $\Xcal$, so (under Assumption~\ref{ass.main}) it follows that $\nabla f$ is bounded in $\ell_2$ (Euclidean) norm by $\kappa_{\nabla f}$ over $\Xcal$.
    \item $\nabla f$ is Lipschitz with constant $L_{\nabla f}$ with respect to $\|\cdot\|_2$ over $\Xcal$.
    \item For all $i \in [m]$, $c_i$ is bounded in absolute value by $\kappa_{c_i}$ over $\Xcal$.
    \item For all $i \in [m]$, $c_i$ is Lipschitz with constant $L_{c_i}$ over $\Xcal$, so (under Assumption~\ref{ass.main}) it follows that $\nabla c_i$ is bounded in $\ell_2$ (Euclidean) norm by $\kappa_{\nabla c_i} := L_{c_i}$ over $\Xcal$.
    \item For all $i \in [m]$, $\nabla c_i$ is Lipschitz with constant $L_{\nabla c_i}$ with respect to $\|\cdot\|_2$ over $\Xcal$.
  \eitemize
\eassumption

Under Assumption \ref{ass.bound}, with $\{\kappa_{c_i}\}_{i\in[m]} \subset \R{}_{>0}$ defined in the assumption, we introduce $\Ccal_{[-\kappa_c,0)} := \{x \in \R{n} : c_i(x) \in [-\kappa_{c_i},0)\ \text{for all}\ i \in [m]\}$ along with the shifted barrier-augmented objective function $\barrier : \Ccal_{[-\kappa_c,0)} \times \R{}_{>0} \to \R{}$ defined by
\bequationNN
  \barrier(x, \mu) = f(x) - \mu \sum_{i\in[m]} \log\(-\frac{c_i(x)}{\kappa_{c_i}}\).
\eequationNN
Important properties of the shifted barrier-augmented objective function are that (a) it is bounded below by $f_{\inf} := \inf_{x \in \Fcal_{<0}} f(x) \in \R{}$ whose existence follows under Assumption~\ref{ass.main} and (b) it shares the same gradient function with respect to its first argument as that of the (unshifted) barrier-augmented objective function; i.e., for all $(x,\mu) \in \Ccal_{[-\kappa_c,0)} \times \R{}_{>0}$, one finds that
\bequation\label{eq.equal_gradients}
  \nabla_x \phi(x,\mu) = \nabla_x \barrier(x,\mu),
\eequation
where $\nabla_x$ denotes the gradient operator with respect to a function's first argument.  These follow since for any $i \in [m]$ and $\delta \in (0,\kappa_{c_i}]$ one has $-\log(\delta/\kappa_{c_i}) \geq 0$ and $-\log(\delta/\kappa_{c_i}) = -\log(\delta) + \log(\kappa_{c_i})$, so $(x,\mu) \in \Ccal_{[-\kappa_{c_i},0)} \times \R{}_{>0}$ implies
\bequation\label{eq.barrier_bounded}
  \phi(x,\mu) + \mu \sum_{i\in[m]} \log(\kappa_{c_i}) = \barrier(x,\mu) \geq f_{\inf}.
\eequation
In addition, the shifted barrier-augmented objective function is monotonically increasing in its second argument in the sense that for any $x \in \Ccal_{[-\kappa_c,0)}$ one finds that $\barrier(x,\mu) < \barrier(x,\bar\mu)$ for any $(\mu,\bar\mu) \in \R{}_{>0} \times \R{}_{>0}$ with $\mu < \bar\mu$.  These properties make the function $\barrier$ useful for our analysis.  Indeed, even though the algorithm is not aware of the constants $\{\kappa_{c_i}\}_{i\in[m]}$ in Assumption~\ref{ass.bound}, we shall show that it follows from the former property above that it computes steps that lead to (expected) decrease in the shifted barrier-augmented objective function $\barrier$.  Also, with the latter property above, decreases in the barrier parameter also lead to decreases in~$\barrier$.

Our first result is useful for both the deterministic and stochastic settings.  It shows that, for any $\theta \in \R{}_{>0}$, the gradient of the shifted barrier-augmented objective function with respect to its first argument is Lipschitz over line segments in $\Ncal(\theta)$ with a Lipschitz constant that depends on $\theta$.  This is a critical property, especially since this function is not globally Lipschitz over $\Fcal_{<0} \supset \Ncal(\theta)$.  The proof follows standard techniques for such a result; we provide it for completeness.

\blemma\label{lem.lipschitz}
  For any $(\mu,\theta,\bar\theta) \in \R{}_{>0} \times \R{}_{>0} \times \R{}_{>0}$ with $\bar\theta \in (0,\theta]$, $x \in \Ncal(\theta)$, $\xbar \in \Ncal(\bar\theta)$, and $\gamma \in [0,1]$, one finds with
  \bequationNN
    L_{\mu,\theta,\bar\theta} := L_{\nabla f} + \frac{\mu}{\theta \bar\theta} \sum_{i=1}^m (L_{c_i} \kappa_{\nabla c_i} + \kappa_{c_i} L_{\nabla c_i}) \in \R{}_{>0}
  \eequationNN
  that one has both
  \begin{align*}
    \|\nabla_x \barrier(x + \gamma (\xbar - x),\mu) - \nabla_x \barrier(x,\mu)\|_2 \leq \gamma L_{\mu,\theta,\bar\theta} \|\xbar - x\|_2 & \\ \text{and}\ \ 
    \barrier(\xbar,\mu) \leq \barrier(x,\mu) + \nabla_x \barrier(x,\mu)^T(\xbar - x) + \thalf L_{\mu,\theta,\bar\theta} \|\xbar - x\|_2^2 &
  \end{align*}
\elemma
\bproof
  For arbitrary $(\mu,\theta,\bar\theta,x,\xbar,\gamma)$ satisfying the conditions of the lemma, it follows from Assumption~\ref{ass.main}, Assumption~\ref{ass.bound}, \eqref{eq.equal_gradients}, and the triangle inequality that
  \begin{align*}
    &\ \|\nabla_x \barrier(x + \gamma(\xbar - x),\mu) - \nabla_x \barrier(x,\mu)\|_2 \\ 
    =&\ \left\|\nabla f(x + \gamma(\xbar - x)) - \nabla f(x) - \mu \( \sum_{i=1}^m \frac{\nabla c_i(x + \gamma(\xbar - x))}{c_i(x + \gamma(\xbar - x))} - \sum_{i=1}^m \frac{\nabla c_i(x)}{c_i(x)} \)\right\|_2 \\
    \leq&\ \gamma L_{\nabla f} \|\xbar - x\|_2 + \mu \sum_{i=1}^m \left\| \frac{ c_i(x + \gamma(\xbar - x))\nabla c_i(x) - c_i(x) \nabla c_i(x + \gamma(\xbar - x))}{c_i(x) c_i(x + \gamma(\xbar - x))} \right\|_2 \\
    \leq&\ \gamma L_{\nabla f} \|\xbar - x\|_2 + \frac{\mu}{\theta \bar\theta} \sum_{i=1}^m \| c_i(x + \gamma(\xbar - x))\nabla c_i(x) - c_i(x) \nabla c_i(x) \\
    &\hspace{100pt} + c_i(x) \nabla c_i(x)  - c_i(x) \nabla c_i(x + \gamma(\xbar - x)) \|_2 \\
    \leq&\ \gamma L_{\nabla f} \|\xbar - x\|_2 + \frac{\mu}{\theta \bar\theta} \sum_{i=1}^m (|c_i(x + \gamma(\xbar - x)) - c_i(x)| \|\nabla c_i(x)\|_2 \\
    &\hspace{100pt} + |c_i(x)|\| \nabla c_i(x) - \nabla c_i(x + \gamma(\xbar - x))\|_2) \\
    \leq&\ \gamma L_{\nabla f} \|\xbar - x\|_2 + \frac{\mu}{\theta \bar\theta} \sum_{i=1}^m (\gamma L_{c_i} \kappa_{\nabla c_i} \|\xbar - x\|_2 + \gamma \kappa_{c_i} L_{\nabla c_i} \|\xbar - x\|_2),
  \end{align*}
  from which the first desired conclusion follows.  Hence, for arbitrary $(\mu,\theta,\bar\theta,x,\xbar)$ satisfying the conditions of the lemma, it follows from the Fundamental Theorem of Calculus and the Cauchy--Schwarz inequality that
  \begin{align*}
    &\ \barrier(\xbar,\mu) - \barrier(x,\mu) \\
    =&\ \int_0^1 \frac{\partial \barrier(x + \gamma(\xbar - x),\mu)}{\partial \gamma} \text{d} \gamma = \int_0^1 \nabla_x \barrier(x + \gamma(\xbar - x),\mu)^T(\xbar - x) \text{d} \gamma \\
    =&\ \nabla_x \barrier(x,\mu)^T(\xbar - x) + \int_0^1 (\nabla_x \barrier(x + \gamma(\xbar - x),\mu) - \nabla_x \barrier(x,\mu))^T(\xbar - x) \text{d} \gamma \\
    \leq&\ \nabla_x \barrier(x,\mu)^T(\xbar - x) + \|\xbar - x\|_2 \int_0^1 \|\nabla_x \barrier(x + \gamma(\xbar - x),\mu) - \nabla_x \barrier(x,\mu))\|_2 \text{d} \gamma \\
    \leq&\ \nabla_x \barrier(x,\mu)^T(\xbar - x) + L_{\mu,\theta,\bar\theta} \|\xbar - x\|_2 \int_0^1 \gamma \text{d} \gamma,
  \end{align*}
  so the desired conclusion follows from the fact that $\int_0^1 \gamma \text{d}\gamma = \thalf$. \qed
\eproof

%************
% Subsection
%************
\subsection{Deterministic Setting}\label{sec.deterministic}

In this subsection, we prove convergence guarantees for Algorithm~\ref{alg.slip} in the deterministic setting, i.e., when $q_k = \nabla f(x_k) - \mu_k \nabla c(x_k) \diag(c(x_k))^{-1} \ones = \nabla_x \barrier(x_k,\mu_k)$ (see \eqref{eq.gradient} and \eqref{eq.equal_gradients}) for all $k \in \N{}$.  For ease of exposition, we refer to this instance of the method as Algorithm~\ref{alg.slip}(D).  In particular, we prove that under a certain rule for choosing the step-size parameters, the algorithm is well defined, generates a sequence of iterates over which a measure of stationarity vanishes, and under a constraint qualification can produce a sequence of Lagrange multipliers such that convergence to a KKT point is guaranteed.  This algorithm for the deterministic setting is not expected to be competitive with state-of-the-art (second-order) interior-point methods for solving problems in many real-world settings.  That said, our analysis for the deterministic setting provides a useful precursor for our subsequent analysis for the stochastic setting for which state-of-the-art methods are not applicable (since they do not possess known convergence guarantees in the stochastic setting that we consider in our analysis in \S\ref{sec.stochastic}).

The class of step-size parameters that we consider for our analysis in this subsection is that satisfying the following parameter rule.  We state the rule, then summarize the reasons that the rule has been designed in this manner.

\begin{parameter}\label{par.deterministic}
  Given $t_\alpha \in (-\infty,0]$, for all $k \in \N{}$, Algorithm~\ref{alg.slip}(D) employs
  \bequationNN
    \alpha_k \gets k^{t_\alpha} \(\frac{\underline\zeta \zeta}{\overline\zeta^2 L_k}\),
  \eequationNN
  where $\displaystyle L_k := L_{\nabla f} + \frac{\mu_k}{\theta_k \theta_{k-1}} \sum_{i=1}^m \(L_{c_i} \kappa_{\nabla c_i} + \kappa_{c_i} L_{\nabla c_i}\) \in \R{}_{>0}$.  In addition, for all $k \in \N{}$, Algorithm~\ref{alg.slip}(D) employs $\gamma_{k,\max} \gets 1$ and $\gamma_k \gets \min_{i\in[m]} \gamma_{k,i}$, where, for all $i \in [m]$,
  \bequationNN
    \gamma_{k,i} \gets \min\left\{1, \frac{-\nabla c_i(x_k)^T d_k + \sqrt{(\nabla c_i(x_k)^T d_k)^2 + 2 L_{\nabla c_i} \|d_k\|_2^2 (-c_i(x_k) - \theta_k)}}{L_{\nabla c_i} \alpha_k \|d_k\|_2^2} \right\};
  \eequationNN
  hence, $\gamma_{k,i} \in (0,1]$ for all $k \in \N{}$ and $i \in [m]$, so $\gamma_k \in (0,1]$ for all $k \in \N{}$.
\end{parameter}

The fundamental aspects of Parameter Rule~\ref{par.deterministic} can be understood as follows.  First, the choice of $\alpha_k$ ensures that this step size is chosen sufficiently small such that it ensures sufficient decrease in the barrier-augmented objective function with each step.  Intuitively, this can be seen in the fact that the step size is proportional to the reciprocal of a Lipschitz-type constant for the gradient of the barrier function (see Lemma~\ref{lem.lipschitz}), as is common in basic descent methods for gradient-based optimization algorithms.  For additional flexibility and since it is required in the stochastic setting, we introduce the factor $k^{t_\alpha}$ such that the step size may diminish faster, namely, when a choice $t_\alpha < 0$ is made.  (That said, our analysis shows that $t_\alpha = 0$ is a valid choice in the deterministic setting.)  Finally, the rule for the choice of $\gamma_k$ is shown in our analysis to guarantee that $[x_k,x_k + \gamma_k \alpha_k d_k] \subseteq \Ncal(\theta_k)$ for all $k \in \N{}$, which is necessary for our analysis to employ Lemma~\ref{lem.lipschitz}.

For our next lemma, we prove this critical property of $\gamma_k$ for all $k \in \N{}$.

\blemma\label{lem.interval}
  Suppose that Assumptions~\ref{ass.main}, \ref{ass.c}, and \ref{ass.bound} hold and that Algorithm~\ref{alg.slip}(D) is run with Parameter Rule~\ref{par.deterministic}.  Then, for all $k \in \N{}$. the choice of $\gamma_k \in (0,1]$ in Parameter Rule~\ref{par.deterministic} ensures that $[x_k,x_k + \gamma_k\alpha_kd_k] \subseteq \Ncal(\theta_k)$, i.e., the condition in Step~\ref{step.gamma} of Algorithm~\ref{alg.slip}(D) holds.  Consequently, under Parameter Rule~\ref{par.deterministic}, Algorithm~\ref{alg.slip}(D) is well defined in the sense that it will generate an infinite sequence of iterates.
\elemma
\bproof
  Consider arbitrary $k \in \N{}$ and $i \in [m]$.  We prove that the choice of $\gamma_{k,i} \in (0,1]$ in Parameter Rule~\ref{par.deterministic} ensures that $c_i(x_k + \gamma \alpha_k d_k) \leq -\theta_k$ for all $\gamma \in [0,\gamma_{k,i}]$, after which the desired conclusion follows directly from the fact that $\gamma_k \gets \min_{i\in[m]} \gamma_{k,i}$ for all $k \in \N{}$.  Under Assumptions~\ref{ass.main} and \ref{ass.bound}, one has for all $\gamma \in [0,1]$ that
  \bequationNN
    c_i(x_k + \gamma \alpha_k d_k) \leq c_i(x_k) + \gamma \alpha_k \nabla c_i(x_k)^T d_k + \thalf L_{\nabla c_i} \gamma^2 \alpha_k^2 \|d_k\|_2^2.
  \eequationNN
  Therefore, $c_i(x_k + \gamma \alpha_k d_k) \leq -\theta_k$ is ensured as long as $\gamma \in [0,1]$ yields
  \bequation\label{eq.c_Lipschitz}
    c_i(x_k) + \bar\gamma \alpha_k \nabla c_i(x_k)^T d_k + \thalf L_{\nabla c_i} \bar\gamma^2 \alpha_k^2 \|d_k\|_2^2 \leq -\theta_k\ \ \text{for all}\ \ \bar\gamma \in [0,\gamma].
  \eequation
  Recalling that Lemma~\ref{lem.in_neighborhoods} ensures $x_k \in \Ncal(\theta_{k-1})$, which means that $c_i(x_k) \leq -\theta_{k-1} < -\theta_k$, it follows that \eqref{eq.c_Lipschitz} holds for all $\gamma \in [0,1]$ when $d_k = 0$, and otherwise (when $d_k \neq 0$) one finds that the left-hand side of \eqref{eq.c_Lipschitz} is a strongly convex quadratic function of $\bar\gamma$.  The second term in the minimum in the definition of $\gamma_{k,i}$ in Parameter Rule~\ref{par.deterministic} is the unique positive root of this quadratic function, where we remark that $-c_i(x_k) - \theta_k > 0$ since, as previously mentioned, $c_i(x_k) < -\theta_k$.
  \qed
\eproof

We now go beyond Lemma~\ref{lem.interval} and prove a lower bound for $\gamma_k$ that will be critical for our ultimate convergence guarantee in this subsection.

\blemma\label{lem.delta}
  Suppose that Assumptions~\ref{ass.main}, \ref{ass.c}, and \ref{ass.bound} hold and that Algorithm~\ref{alg.slip}(D) is run with Parameter Rule~\ref{par.deterministic}.  In addition, define
  \bequationNN
    \beta := \overline\zeta \(\kappa_{\nabla f} + \frac{\mu_1}{\theta_0} \sum_{j \in [m]} \kappa_{\nabla c_j}\) \in \R{}_{>0},
  \eequationNN
  and for all $i \in [m]$ define
  \bequationNN
    \gamma_{k,i,\min} := \min\left\{1, \frac{-\kappa_{\nabla c_i} + \sqrt{\kappa_{\nabla c_i}^2 + 2 L_{\nabla c_i} (\eta \mu_k - \theta_k)}}{\alpha_k \beta L_{\nabla c_i}}, \frac{\underline\eta - \theta_k}{\alpha_k \beta L_{\nabla c_i}} \right\} \in (0,1].
  \eequationNN
  Then, $\gamma_k \geq \gamma_{k,\min} := \min_{i\in[m]} \gamma_{k,i,\min}$ for all $k \in \N{}$.
\elemma
\bproof
  Consider arbitrary $k \in \N{}$ and $i \in [m]$.  Our aim is to prove that $\gamma_{k,i} \geq \gamma_{k,i,\min}$, from which the desired conclusion follows due to the fact that $\gamma_k \gets \min_{i\in[m]} \gamma_{k,i}$ under Parameter Rule~\ref{par.deterministic}.  Toward this end, observe that if Parameter Rule~\ref{par.deterministic} yields $\gamma_{k,i} = 1$, then, clearly, $\gamma_{k,i} \geq \gamma_{k,i,\min} \in (0,1]$ and there is nothing left to prove.  Therefore, we may proceed under the assumption that Parameter Rule~\ref{par.deterministic} yields $\gamma_{k,i} \in (0,1)$.  It follows under this condition that
  \bequation\label{eq.qi}
    \gamma_{k,i} = \frac{\frac{-\nabla c_i(x_k)^Td_k}{\|d_k\|_2} + \sqrt{\(\frac{\nabla c_i(x_k)^Td_k}{\|d_k\|_2}\)^2 + 2 L_{\nabla c_i} (-c_i(x_k) - \theta_k)}}{\alpha_k \|d_k\|_2 L_{\nabla c_i}}.
  \eequation
  On the other hand, by Assumptions~\ref{ass.main}, \ref{ass.c}, and \ref{ass.bound}, the fact that the matrix norm $\|\cdot\|_2$ is submultiplicative, the triangle inequality, \eqref{eq.d_condition_3}, the fact that $x_k \in \Ncal(\theta_{k-1})$, the fact that $\{\mu_k/\theta_{k-1}\}$ is a constant sequence, and the fact that for the orthogonal projection matrix $P$ one has $\|P\|_2 \leq 1$, one finds that
  \bequation\label{eq.beta}
    \baligned
      \|d_k\|_2
        &\leq \overline\zeta \|P(\nabla f(x_k) - \mu_k \nabla c(x_k) \diag(c(x_k))^{-1} \ones)\|_2 \\
        &\leq \overline\zeta \(\|\nabla f(x_k)\|_2 + \mu_k \sum_{j\in[m]} \left\| c_j(x_k)^{-1} \nabla c_j(x_k) \right\|_2 \) \leq \beta.
    \ealigned
  \eequation
  Let us now proceed by considering two cases.
  \benumerate
    \item Suppose $c_i(x_k) \leq -\eta \mu_k$. (The $i$th constraint is not nearly active.) Observe that Assumption~\ref{ass.bound} and the Cauchy-Schwarz inequality together imply that $\nabla c_i(x_k)^Td_k \leq \kappa_{\nabla c_i} \|d_k\|_2$.  At the same time, observe that the numerator of the right-hand side of \eqref{eq.qi} can be viewed as the value of the function $h : (-\infty,\kappa_{\nabla c_i}] \to \R{}$ defined by $h(a) = -a + \sqrt{a^2 + b}$ at $a = \frac{\nabla c_i(x_k)^Td_k}{\|d_k\|_2}$ for $b \in \R{}_{>0}$, where positivity of $b$ follows since $x_k \in \Ncal(\theta_{k-1})$ implies $-c_i(x_k) - \theta_k \geq \theta_{k-1} - \theta_k > 0$.  One finds that $h$ is a monotonically decreasing function of~$a$ over $(-\infty,\kappa_{\nabla c_i}]$.  Hence, from \eqref{eq.qi}, \eqref{eq.beta}, and the condition of this case, one has
    \bequation\label{eq.qi_2}
      \baligned
        \gamma_{k,i} &\geq \frac{-\kappa_{\nabla c_i} + \sqrt{\kappa_{\nabla c_i}^2 + 2 L_{\nabla c_i} (-c_i(x_k) - \theta_k)}}{\alpha_k \beta L_{\nabla c_i}} \\
        &\geq \frac{-\kappa_{\nabla c_i} + \sqrt{\kappa_{\nabla c_i}^2 + 2 L_{\nabla c_i} (\eta \mu_k - \theta_k)}}{\alpha_k \beta L_{\nabla c_i}}.
      \ealigned
    \eequation
    \item Suppose $c_i(x_k) > -\eta \mu_k$. (The $i$th constraint is nearly active.)  Under Assumption~\ref{ass.c}, it follows from \eqref{eq.d_condition_5} that $\nabla c_i(x_k)^Td_k \leq - \thalf \underline\eta \|d_k\|_2$.  Consequently, by $\theta_k > 0$, the fact that $x_k \in \Ncal(\theta_{k-1})$ implies $-c_i(x_k) - \theta_k > 0$, and the fact that the conditions of the lemma ensure $\underline\eta - \theta_k > 0$, it follows that
    \begin{align*}
      &&\!\!\!\!\!\! -\frac{\underline\eta}{2} &\geq \frac{\nabla c_i(x_k)^Td_k}{\|d_k\|_2} \\
      \implies &&\!\!\!\!\!\! L_{\nabla c_i} \(\frac{-c_i(x_k) - \theta_k}{\underline\eta - \theta_k}\) - \frac{\underline\eta - \theta_k}{2} &\geq \frac{\nabla c_i(x_k)^Td_k}{\|d_k\|_2} \\
      \implies &&\!\!\!\!\!\! 2L_{\nabla c_i} (-c_i(x_k) - \theta_k) - (\underline\eta - \theta_k)^2 &\geq 2(\underline\eta - \theta_k) \frac{\nabla c_i(x_k)^Td_k}{\|d_k\|_2} \\
      \implies &&\!\!\!\!\!\! \(\frac{\nabla c_i(x_k)^Td_k}{\|d_k\|_2}\)^2 + 2L_{\nabla c_i} (-c_i(x_k) - \theta_k) &\geq \(\frac{\nabla c_i(x_k)^Td_k}{\|d_k\|_2} + (\underline\eta - \theta_k) \)^2 \\
      \implies &&\!\!\!\!\!\! \sqrt{ \(\frac{\nabla c_i(x_k)^Td_k}{\|d_k\|_2}\)^2 + 2L_{\nabla c_i} (-c_i(x_k) - \theta_k) } &\geq \frac{\nabla c_i(x_k)^Td_k}{\|d_k\|_2} + (\underline\eta - \theta_k),
    \end{align*}
    which with \eqref{eq.qi} and \eqref{eq.beta} implies
    \bequationNN
      \gamma_{k,i} \geq \frac{\underline\eta - \theta_k}{\alpha_k \beta L_{\nabla c_i}}.
    \eequationNN
  \eenumerate
  Combining the results from the two cases, one can conclude that $\gamma_{k,i} \geq \gamma_{k,i,\min}$ for all $i \in [m]$, which, as previously mentioned, completes the proof. \qed
\eproof

For our next lemma, we prove the aforementioned sufficient decrease condition that is guaranteed by the choices in Parameter Rule~\ref{par.deterministic}.  In particular, one finds in the lemma that the amount of decrease is at least proportional to a squared norm of the \emph{projected} gradient of the barrier-augmented objective function, which makes sense since, under Lemma~\ref{lem.in_neighborhoods}, $x_k \in \Ecal$ and $d_k \in \Null(A)$ for all $k \in \N{}$.

\blemma\label{lem.sufficient_decrease}
  Suppose that Assumptions~\ref{ass.main}, \ref{ass.c}, and \ref{ass.bound} hold and that Algorithm~\ref{alg.slip}(D) is run with Parameter Rule~\ref{par.deterministic}.  Then, for all $k \in \N{}$, one finds
  \bequationNN
    \barrier(x_{k+1}, \mu_{k+1}) -  \barrier(x_k, \mu_k) \leq - \thalf \underline\zeta \zeta \gamma_k \alpha_k \| P \nabla_x \phi(x_k, \mu_k) \|_2^2.
  \eequationNN
\elemma
\bproof
  Consider arbitrary $k \in \N{}$.  By Lemmas~\ref{lem.lipschitz}--\ref{lem.interval}, \eqref{eq.equal_gradients}, \eqref{eq.d_conditions}, $d_k = Pd_k$, and $P = P^T$, it follows with $L_k \in \R{}_{>0}$ defined as in Parameter Rule~\ref{par.deterministic} that
  \begin{align}
    \barrier(x_{k+1}, \mu_k) -  \barrier(x_k, \mu_k)
    \leq&\ \nabla_x \barrier(x_k, \mu_k)^T (x_{k+1} - x_k) + \thalf L_k \|x_{k+1} - x_k\|_2^2 \nonumber \\
    =&\ \gamma_k \alpha_k q_k^T d_k + \thalf \gamma_k^2 \alpha_k^2 L_k \|d_k\|_2^2 \nonumber \\
    =&\ \gamma_k \alpha_k (P q_k)^Td_k + \thalf \gamma_k^2 \alpha_k^2 L_k \|d_k\|_2^2 \nonumber \\
    \leq&\ -\gamma_k \alpha_k \zeta \|Pq_k\|_2 \|d_k\|_2 + \thalf \gamma_k^2 \alpha_k^2 L_k \|d_k\|_2^2 \nonumber \\
    \leq&\ -\gamma_k \alpha_k (\underline\zeta \zeta - \thalf \overline\zeta^2 \gamma_k \alpha_k L_k) \| Pq_k \|_2^2. \label{eq.suff1}
  \end{align}
  Due to Parameter Rule~\ref{par.deterministic}, one finds that $\gamma_k \in (0,1]$, $k^{t_\alpha} \in (0,1]$, and
  \bequation\label{eq.suff2}
    \underline\zeta \zeta - \thalf \overline\zeta^2 \gamma_k \alpha_k L_k \geq \underline\zeta \zeta - \thalf \overline\zeta^2 \alpha_k L_k \geq \underline\zeta \zeta - \thalf k^{t_\alpha} \underline\zeta \zeta \geq \thalf \underline\zeta \zeta.
  \eequation
  Combining \eqref{eq.suff1} and \eqref{eq.suff2}, and since $\mu_{k+1} \leq \mu_k$ implies $\barrier (x_{k+1}, \mu_{k+1}) \leq \barrier (x_{k+1}, \mu_k)$ (see the discussion following~\eqref{eq.barrier_bounded}), the desired conclusion follows. \qed
\eproof

We are now prepared to prove our guarantee for the deterministic setting.  For the sake of notational clarity, the following theorem introduces the parameters $t_\mu$ and $t_\theta$ that dictate the rates of decrease of $\{\mu_k\}$ and $\{\theta_k\}$, respectively.  However, the theorem requires $t_\mu = t_\theta$, and indeed one can see in the details of the proof that both $t_\mu \leq t_\theta$ and $t_\mu \geq t_\theta$ are required to prove the theorem.

\btheorem\label{th.deterministic}
  Suppose that Assumptions~\ref{ass.main}, \ref{ass.c}, and \ref{ass.bound} hold and that Algorithm~\ref{alg.slip}(D) is run with Parameter Rule~\ref{par.deterministic}.  Suppose, in addition, that for some $(t_\mu,t_\theta,t_\alpha) \in (-\infty,0) \times (-\infty,0) \times (-\infty,0]$ with $t_\mu = t_\theta$ and $t_\mu + t_\alpha \in [-1,0)$ and for some $\mu_1 \in \R{}_{>0}$ and $\theta_0 \in \R{}_{>0}$ the algorithm employs $\mu_k = \mu_1 k^{t_\mu}$ and $\theta_{k-1} = \theta_0 k^{t_\theta}$ for all $k \in \N{}$.  Then,
  \bequation\label{eq.liminf}
    \sum_{k=1}^\infty \gamma_k \alpha_k = \infty\ \ \text{and}\ \ \liminf_{k\to\infty} \|P \nabla_x \phi(x_k,\mu_k)\|_2^2 = 0.
  \eequation
  If, in addition, there exists $\Kcal \subseteq \N{}$ with $|\Kcal| = \infty$ such that
  \bitemize
    \renewcommand\labelitemi{$\bullet$}
    \item $\{P\nabla_x \phi(x_k,\mu_k)\}_{k\in\Kcal} \to 0$,
    \item $\{x_k\}_{k\in\Kcal} \to \xbar$ for some $\xbar \in \Fcal$, and
    \item at $\xbar$ the linear independence constraint qualification (LICQ) holds with respect to problem~\eqref{prob.opt} in the sense that with $\bar\Acal := \{i \in [m] : c_i(\xbar) = 0\}$ the columns of $A^T$ combined with the vectors in $\{\nabla c_i(\xbar)\}_{i\in\bar\Acal}$ form a linearly independent set,
  \eitemize
  then $\xbar$ is a KKT point for problem~\eqref{prob.opt} in the sense that there exists a pair of Lagrange multipliers $(\ybar,\zbar) \in \R{l} \times \R{m}$ such that $(\xbar,\ybar,\zbar)$ satisfies \eqref{eq.KKT}.
\etheorem
\bproof
  It follows from \eqref{eq.barrier_bounded} that $\barrier$ is bounded below by $f_{\inf}$ over $\Xcal \times \R{}_{\geq0}$.  Then, by summing the expression in Lemma~\ref{lem.sufficient_decrease} over $k \in \N{}$ and \eqref{eq.equal_gradients}, one finds that
  \begin{align*}
    \infty > \barrier(x_1,\mu_1) - f_{\inf} &\geq \sum_{k=1}^\infty (\barrier(x_k,\mu_k) - \barrier(x_{k+1},\mu_{k+1})) \\
    &\geq \sum_{k=1}^\infty \thalf \underline\zeta \zeta \gamma_k \alpha_k \|P \nabla_x \phi(x_k,\mu_k)\|_2^2.
  \end{align*}
  To complete the proof of \eqref{eq.liminf}, let us now show that $\{\gamma_k \alpha_k\}$ is unsummable.  To begin, first observe that the lower bound on $\{\gamma_k\}$ stated in Lemma~\ref{lem.delta} holds, i.e., there exist $\kappa \in \R{}_{>0}$, $\underline{L} \in \R{}_{>0}$, and $L \in \R{}_{>0}$ such that, for all $k \in \N{}$, one has
  \bequationNN
    \gamma_{k,\min} \geq \min\left\{1, \frac{-\kappa + \sqrt{\kappa^2 + 2 \underline{L} (\eta \mu_k - \theta_k)}}{\alpha_k \beta L}, \frac{\underline\eta - \theta_k}{\alpha_k \beta L} \right\}.
  \eequationNN
  Consequently, it holds for all $k \in \N{}$ that
  \begin{align}
    \gamma_k \alpha_k \geq \gamma_{k,\min} \alpha_k
    &\geq \min\left\{ \alpha_k, \frac{-\kappa + \sqrt{\kappa^2 + 2 \underline{L} (\eta \mu_k - \theta_k)}}{\beta L}, \frac{\underline\eta - \theta_k}{\beta L} \right\} \nonumber \\
    &=: \min\{ \alpha_k, \rho_k, \upsilon_k \}. \label{def:betak-alphak}
  \end{align}
  Let us now consider the sequences $\{\alpha_k\}$, $\{\rho_k\}$, and $\{\upsilon_k\}$ in turn.  Our aim is to show that each of these sequences is unsummable, which will complete the proof of~\eqref{eq.liminf}.  With respect to the step-size sequence $\{\alpha_k\}$, one finds by Parameter Rule~\ref{par.deterministic}, $\upsilon := \sum_{i=1}^m \(L_{c_i} \kappa_{\nabla c_i} + \kappa_{c_i} L_{\nabla c_i}\)$, and the conditions of the theorem that
  \bequation\label{eq.alpha_lower}
    \alpha_k = \frac{k^{t_\alpha} \underline\zeta \zeta \overline\zeta^{-2}}{L_{\nabla f} + \mu_1 \theta_0^{-2} k^{t_\mu} k^{-t_\theta} (k + 1)^{-k_\theta} \upsilon}.
  \eequation
  Observe that $t_\mu \in (-\infty,0)$, $t_\theta \in (-\infty,0)$, and $t_\mu = t_\theta$ show that $k^{t_\mu} k^{-t_\theta} (k + 1)^{-t_\theta} \leq k^{t_\mu-t_\theta} (2k)^{-t_\theta} = 2^{-t_\theta} k^{t_\mu - 2t_\theta} = 2^{-t_\mu} k^{-t_\mu}$.  With \eqref{eq.alpha_lower} and $t_\mu \in (-\infty,0)$, this shows
  \bequationNN
    \alpha_k \geq \frac{k^{t_\alpha} \underline\zeta \zeta \overline\zeta^{-2}}{L_{\nabla f} + \mu_1 \theta_0^{-2} 2^{-t_\mu} \upsilon k^{-t_\mu}} = \frac{k^{t_\mu + t_\alpha} \underline\zeta \zeta \overline\zeta^{-2}}{L_{\nabla f} k^{t_\mu} + \mu_1 \theta_0^{-2} 2^{-t_\mu} \upsilon} \geq \frac{k^{t_\mu + t_\alpha} \underline\zeta \zeta \overline\zeta^{-2}}{L_{\nabla f} + \mu_1 \theta_0^{-2} 2^{-t_\mu} \upsilon},
  \eequationNN
  which since $t_\mu + t_\alpha \in [-1,0)$ implies that $\{\alpha_k\}$ is unsummable, as desired.  Now, with respect to $\{\rho_k\}$, one finds under the conditions of the theorem that
  \bequation\label{eq.rho_lower}
    \rho_k = \frac{-\kappa + \sqrt{\kappa^2 + 2 \underline{L} (\eta \mu_1 k^{t_\mu} - \theta_0 (k + 1)^{t_\theta})}}{\beta L}.
  \eequation
  Observe that $t_\mu \in (-\infty,0)$ and $t_\mu = t_\theta$ show that
  \bequationNN
    \eta \mu_1 k^{t_\mu} - \theta_0 (k + 1)^{t_\theta} \geq \eta \mu_1 k^{t_\mu} - \theta_0 k^{t_\theta} = (\eta \mu_1 - \theta_0) k^{t_\mu}.
  \eequationNN
  In fact, $t_\mu \in [-1,0)$, so the bound above along with \eqref{eq.rho_lower} shows that $\{\rho_k\}$ is unsummable, as desired.  Finally, $\{\upsilon_k\}$ can be seen to be unsummable since $\upsilon_k = \underline\eta - \theta_k > \underline\eta - \theta_0$ for all $k \in \N{}$.  Having reached the desired conclusion that $\{\alpha_k\}$, $\{\rho_k\}$, and $\{\upsilon_k\}$ are unsummable, it follows through \eqref{def:betak-alphak} that $\{\gamma_k \alpha_k\}$ is unsummable, which, as previously mentioned, completes the proof of \eqref{eq.liminf}.

  Now suppose that there exists an infinite-cardinality set $\Kcal \subseteq \N{}$ as described in the theorem.  By construction of the algorithm, it follows that $Ax_k = b$ and $c(x_k) < 0$ for all $k \in \N{}$.  Now, for all $k \in \Kcal$, define the Lagrange multiplier estimates
  \bequation\label{eq.multipliers}
    z_k := -\mu_k \diag(c(x_k))^{-1} \ones \ \ \text{and}\ \ y_k := -(AA^T)^{-1} A (\nabla f(x_k) + \nabla c(x_k)z_k).
  \eequation
  Since $\mu_k > 0$ and $c(x_k) < 0$ for all $k \in \N{}$, it follows that $z_k \geq 0$ for all $k \in \N{}$.  In addition, the fact that $\{P \nabla \phi(x_k,\mu_k)\}_{k\in\Kcal} \to 0$ shows that
  \begin{align}
     &\ \{\|\nabla f(x_k) + A^Ty_k + \nabla c(x_k)z_k\|_2\}_{k \in \Kcal} \nonumber \\
    =&\ \{\|(\nabla f(x_k) + \nabla c(x_k)z_k) - A^T (AA^T)^{-1} A (\nabla f(x_k) + \nabla c(x_k)z_k)\|_2\}_{k \in \Kcal} \nonumber \\
    =&\ \{\|P \nabla \phi(x_k,\mu_k)\|_2\}_{k\in\Kcal} \to 0. \label{eq:xk-to-xbar}
  \end{align}
  Consider $\bar\Acal = \{i \in [m] : c_i(\xbar) = 0\}$.  Since $c(x_k) \circ z_k = - \mu_k \ones$ for all $k \in \N{}$ and $\{\mu_k\} \searrow 0$, it follows that $\{[z_k]_i\}_{k \in \Kcal} \to 0$ for all $i \not\in \bar\Acal$.  Combining this with \eqref{eq:xk-to-xbar} and using the fact that $\{\nabla c_i\}_{k \in \Kcal}$ is bounded for all $i \in [m]$ under Assumption~\ref{ass.bound},
  \bequationNN
    \bigg\{ \bigg\| \nabla f(x_k) + A^Ty_k + \sum_{i\in\bar\Acal} \nabla c_i(x_k) [z_k]_i \bigg\|_2 \bigg\}_{k\in\Kcal} \to 0.
  \eequationNN
  Under the LICQ, it follows from this limit that $\{(y_k,z_k)\}_{k \in \Kcal}$ converges to some pair $(\ybar,\zbar)$ such that $(\xbar,\ybar,\zbar)$ satisfies \eqref{eq.KKT}, as desired.
  \qed
\eproof

%************
% Subsection
%************
\subsection{Stochastic Setting}\label{sec.stochastic}

We now prove convergence guarantees for Algorithm~\ref{alg.slip} in a stochastic setting, i.e., when $q_k = g_k - \mu_k \nabla c(x_k) \diag(c(x_k))^{-1} \ones$ (see \eqref{eq.gradient}) for all $k \in \N{}$ where $g_k$ is a stochastic estimate of $\nabla f(x_k)$.  Formally, let us consider the probability space $(\Omega, \Gcal, \P)$, where $\Omega$ captures all outcomes of a run of Algorithm~\ref{alg.slip}.  As mentioned shortly, we make assumptions that guarantee that each iteration is well defined, which implies that each such outcome corresponds to an infinite sequence of iterates.  In this manner, each realization of a run of the algorithm can be associated with $\omega \in \Omega$, an infinite-dimensional tuple whose $k$th element determines the stochastic gradient-of-the-objective estimate in iteration $k \in \N{}$.  The stochastic process defined by the algorithm can thus be expressed as
\bequationNN
  \{(X_k(\omega), G_k(\omega), Q_k(\omega), D_k(\omega), \Acal_k(\omega), \Gamma_k(\omega))\},
\eequationNN
where, for all $k \in \N{}$, the random variables are the iterate $X_k(\omega)$, stochastic objective-gradient estimator $G_k(\omega)$, stochastic gradient estimator
\bequation\label{eq.Q}
  Q_k(\omega) := G_k(\omega) - \mu_k \nabla c(X_k(\omega)) \diag(c(X_k(\omega)))^{-1} \ones,
\eequation
direction $D_k(\omega)$, step size~$\Acal_k(\omega)$, and neighborhood-enforcement parameter $\Gamma_k(\omega)$.  Observe that the statement of Algorithm~\ref{alg.slip} on page~\pageref{alg.slip} instantiates a particular realization of a run, expressed as $\{(x_k, g_k, q_k, d_k, \alpha_k, \gamma_k)\}$.

Given the initial conditions of the algorithm (including for simplicity that $X_1(\omega) = x_1$ for all $\omega \in \Omega$), let $\Gcal_1$ denote the sub-$\sigma$-algebra of $\Gcal$ corresponding to the initial conditions and, for all $k \in \N{} \setminus \{1\}$, let $\Gcal_k$ denote the sub-$\sigma$-algebra of~$\Gcal$ corresponding to the initial conditions and $\{G_1(\omega),\dots,G_{k-1}(\omega)\}$.  In this manner, the sequence $\{\Gcal_k\}$ is a filtration.  For our analysis, with respect to this filtration, we make the following assumption.  For the sake of brevity, here and for the remainder of our analysis, we drop from our notation the dependence of the stochastic process on an element $\omega$ of $\Omega$, since this dependence remains obvious.

\bassumption\label{ass.stochastic}
  For all $k \in \N{}$, it holds that $\E[G_k | \Gcal_k] = \nabla f(X_k)$. In addition, there exists $\sigma \in \R{}_{\geq 0}$ such that, for all $k \in \N{}$, one has $\|P(G_k - \nabla f(X_k))\|_2 \leq \sigma$.
\eassumption

Assumption~\ref{ass.stochastic} amounts to the stochastic-gradient estimators being unbiased with bounded error.  This is consistent with the stochastic setting in \cite{CKRW23}.  Looser assumptions are possible for stochastic-gradient-based methods in the unconstrained setting, but for the context of stochastic-gradient-based interior-point methods we know of no analysis that can avoid a bounded-error assumption.

We also carry forward our prior assumptions, namely, Assumptions~\ref{ass.main}, \ref{ass.c}, and~\ref{ass.bound}.  However, we strengthen Assumption~\ref{ass.c} somewhat to impose additional structure on the manner in which the search direction is defined.  In particular, for all $k \in \N{}$, we define the search direction $D_k$ through the linear system
\bequation\label{eq.linsys}
  \bbmatrix H_k & A^T \\ A & 0 \ebmatrix \bbmatrix D_k \\ Y_k \ebmatrix = - \bbmatrix Q_k \\ 0 \ebmatrix,
\eequation
where the sequence $\{H_k\}$ satisfies the following assumption.

\bassumption\label{ass.H}
  For all $k \in \N{}$, the matrix $H_k \in \mathbb{S}_{\succeq 0}^n$ is $\Gcal_k$-measurable and $D_k$ satisfies~\eqref{eq.linsys}.  In addition, there exist constants $(\underline\lambda, \overline\lambda) \in (0,1] \times \R{}_{>0}$ with $\underline\lambda \leq \overline\lambda$ such that, over any run of Algorithm~\ref{alg.slip}, one finds $u^TH_ku/\|u\|^2 \in [\underline\lambda, \overline\lambda]$ for all $u \in \Null(A)$.
\eassumption

All combined, Assumptions~\ref{ass.c} and \ref{ass.H} amount to the assumption that $D_k$ is computed by \eqref{eq.linsys} with a choice of $H_k$ such that \eqref{eq.d_condition_5} holds.  Here, we are remarking on the fact that, under Assumption~\ref{ass.H}, the search direction satisfies the null-space and angle conditions imposed by the combination of \eqref{eq.d_condition_1}--\eqref{eq.d_condition_4}.  Indeed, one finds
\bequation\label{eq.DZQ}
  D_k = - Z(Z^TH_kZ)^{-1} Z^TQ_k,
\eequation
where $Z \in \R{n \times (n-l)}$ is an orthogonal matrix whose columns form an orthonormal basis for $\Null(A)$, through which one can in turn show that $AD_k = 0$ along with
\bequation\label{eq.d_linsys_bounds}
  \|D_k\|_2 \in [\overline\lambda^{-1} \|PQ_k\|_2,\underline\lambda^{-1} \|PQ_k\|_2]\ \ \text{and}\ \ -(PQ_k)^TD_k \geq \underline\lambda \overline\lambda^{-1} \|PQ_k\|_2 \|D_k\|_2.
\eequation
Hence, we employ $(\underline\lambda,\overline\lambda)$ and \eqref{eq.d_linsys_bounds} in place of $(\underline\zeta,\overline\zeta,\zeta)$ and \eqref{eq.d_condition_2}--\eqref{eq.d_condition_4} in this section.  For our analysis, it is worthwhile to emphasize that the tuple of parameters $(\eta,\bar\eta,\underline\lambda,\overline\lambda)$ is presumed to be determined uniquely for any given instance of problem~\eqref{prob.opt}.  Therefore, the parameters in Assumption~\ref{ass.c} (i.e., in \eqref{eq.d_conditions}) and (as explicitly stated) in Assumption~\ref{ass.H} are presumed to be uniform over all possible runs of the algorithm.  Similarly, for Assumption~\ref{ass.bound} in this section, the set $\Xcal$ and stated constants are assumed to be uniform over all possible realizations of $\{X_k\}$.

Our convergence guarantees for this stochastic setting rely on Parameter Rule~\ref{par.stochastic}, stated below, that we employ for choosing step-size- and neighborhood-enforcement parameter sequences, namely, $\{\alpha_k\}$ and $\{(\gamma_{k,\min},\gamma_k,\gamma_{k,\max})\}$.  Like in the deterministic setting, in this stochastic setting, the step sizes depend on a parameter $t_\alpha \in (-\infty,0)$, with one difference being that the choice $t_\alpha = 0$ is not allowed for the stochastic setting.  However, unlike in the deterministic setting, in this stochastic setting we prescribe a particular dependence of the neighborhood-enforcement values $\{(\gamma_{k,\min},\gamma_k,\gamma_{k,\max})\}$ on the barrier-parameter sequence $\{\mu_k\}$ and neighborhood-parameter sequence $\{\theta_k\}$, which in turn are determined by parameters $(t_\mu,t_\theta) \in (-\infty,0) \times (-\infty,0)$.  Whereas in the deterministic setting our ultimate requirements for the tuple $(t_\mu, t_\theta, t_\alpha)$ (specified between Parameter Rule~\ref{par.deterministic} and Theorem~\ref{th.deterministic}) were that $t_\mu = t_\theta < 0$, $t_\alpha \leq 0$, and $t_\mu + t_\alpha \in [-1,0)$, the requirements for these parameters is stricter in this stochastic setting.  For the sake of simplicity, rather than introduce separate parameters $t_\mu$ and $t_\theta$ for $\{\mu_k\}$ and~$\{\theta_k\}$, respectively, in Parameter Rule~\ref{par.stochastic} we introduce a single parameter $t$.  This is reasonable since, even in the deterministic setting, our analysis requires $t_\mu = t_\theta$.  It is worthwhile to emphasize upfront that the restrictions of the rule are satisfiable; e.g., one may consider $t = -0.99$ and $t_\alpha = -0.01$ as one acceptable setting.

\begin{parameter}\label{par.stochastic}
  Suppose the following are given:
  \bitemize
    \renewcommand\labelitemi{$\bullet$}
    \item $(t, t_\alpha) \in (-\infty, 0)^2$ such that $t + t_\alpha \in [-1,0)$, $2t + t_\alpha < -1$, and $t + 2t_\alpha < -1$;
    \item $\mu_1 \in \R{}_{>0}$ and $\theta_0 \in \R{}_{>0}$; and
    \item $\gamma_\buff \in \R{}_{\geq 0}$ and $\{\gamma_{k,\buff}\} \subset \R{}_{\geq 0}$ such that $\gamma_{k,\buff} \leq \gamma_\buff k^t$ for all $k \in \N{}$.
  \eitemize
  Algorithm~\ref{alg.slip} employs, for all $k \in \N{}$, the following:
  \bitemize
    \renewcommand\labelitemi{$\bullet$}
    \item $\mu_k = \mu_1 k^t$ and $\theta_{k-1} = \theta_0 k^t$;
    \item $\displaystyle \alpha_k \gets k^{t_\alpha} \(\frac{\underline\lambda^2}{\overline\lambda L_k}\)$, where $L_k \in \R{}_{>0}$ is defined as in Parameter Rule~\ref{par.deterministic}; and
    \item $\gamma_k \gets \min_{i \in [m]} \gamma_{k,i}$ where
    \begin{align*}
      \gamma_{k,i} &\gets \min\{\tilde\gamma_{k,i}, \gamma_{k,\max}\}\ \ \text{for all $i \in [m]$}, \\
      \tilde\gamma_{k,i} &\gets \frac{-\nabla c_i(x_k)^T d_k + \sqrt{(\nabla c_i(x_k)^T d_k)^2 + 2 L_{\nabla c_i} \|d_k\|_2^2 (-c_i(x_k) - \theta_k)}}{L_{\nabla c_i} \alpha_k \|d_k\|_2^2} \\
      &\qquad \text{for all $i \in [m]$}, \\
      \beta^\sigma &\gets \underline\lambda^{-1} \(\kappa_{\nabla f} + \sigma + \frac{\mu_1}{\theta_0} \sum_{j \in [m]} \kappa_{\nabla c_j}\), \\
      \gamma_{k,i,\min} &\gets \min\left\{1, \frac{-\kappa_{\nabla c_i} + \sqrt{\kappa_{\nabla c_i}^2 + 2 L_{\nabla c_i} (\eta \mu_k - \theta_k)}}{\alpha_k \beta^\sigma L_{\nabla c_i}}, \frac{\underline\eta - \theta_k}{\alpha_k \beta^\sigma L_{\nabla c_i}} \right\} \\
      &\qquad \text{for all $i \in [m]$}, \\
      \gamma_{k,\min} &\gets \min_{i \in [m]} \gamma_{k,i,\min},\ \ \text{and} \\
      \gamma_{k,\max} &\gets \min \{1, \gamma_{k,\min} + \gamma_{k,\buff} \}.
    \end{align*}
  \eitemize
\end{parameter}

The strategy for selecting $\{\alpha_k\}$ in Parameter Rule~\ref{par.stochastic} is consistent with that in Parameter Rule~\ref{par.deterministic}.  Observe that since $\{\mu_k\}$ and $\{\theta_k\}$ are set in advance of any run of the algorithm, it follows under Parameter Rule~\ref{par.stochastic} that $\{L_k\}$ and so $\{\alpha_k\}$ are determined in advance of a run of the algorithm as well.  Thus, our analysis can consider $\{\Acal_k\} = \{\alpha_k\}$, where $\{\alpha_k\}$ is a prescribed sequence that is uniform over all runs of Algorithm~\ref{alg.slip} employed for any given instance of problem~\eqref{prob.opt}.  The choice of $\{(\gamma_{k,\min},\gamma_k,\gamma_{k,\max})\}$ in Parameter Rule~\ref{par.stochastic}, on the other hand, is more stringent than that in Parameter Rule~\ref{par.deterministic}.  Observe that Parameter Rule~\ref{par.deterministic} amounts to choosing $\gamma_k \gets \min_{i\in[m]} \gamma_{k,i}$ for all $k \in \N{}$, where $\gamma_{k,i} \gets \min\{\tilde\gamma_{k,i}, 1\}$ for all $(k,i) \in \N{} \times [m]$.  The choice in Parameter Rule~\ref{par.stochastic} is similar, but with $\gamma_{k,i} \gets \min\{\tilde\gamma_{k,i}, \gamma_{k,\max}\}$ for all $(k,i) \in \N{} \times [m]$, where for any given $k \in \N{}$ the upper limit $\gamma_{k,\max} \in (0,1]$ might be less than 1.  The particular strategy for setting $\gamma_{k,\max}$ for all $k \in \N{}$ in Parameter Rule~\ref{par.stochastic} involves first setting a lower value $\gamma_{k,\min} \in \R{}_{>0}$ for each $k \in \N{}$ \emph{that is defined only by values that are set in advance of a run of the algorithm}.  This fact is critical, since it ensures that $\{(\gamma_{k,\min},\gamma_{k,\max})\}$ is prescribed and uniform over all runs of Algorithm~\ref{alg.slip} employed for any given instance of problem~\eqref{prob.opt}.  For each $(k,i) \in \N{} \times [m]$, the values $\tilde\gamma_{k,i}$ and $\gamma_k$ may differ between runs, but importantly the sequence $\{(\gamma_{k,\min},\gamma_{k,\max})\}$ does not.

Many of our prior results from the deterministic setting carry over here to our present analysis for the stochastic setting.  Firstly, Lemma~\ref{lem.in_neighborhoods} holds, as does Lemma~\ref{lem.lipschitz} since it has been proved irrespective of any particular algorithm.  Secondly, following the arguments in Lemmas~\ref{lem.interval} and \ref{lem.delta}, it follows that Parameter Rule~\ref{par.stochastic} ensures that $\gamma_{k,\min} \leq \gamma_k \leq \gamma_{k,\max}$ for all $k \in \N{}$ in any given run of the algorithm.  We formalize this in Lemma~\ref{lem.delta_stochastic} below, the proof of which is omitted since it would follow the same lines of arguments as in the proofs of Lemmas~\ref{lem.interval} and \ref{lem.delta}.  The only significant difference that would be in the proof is that the bound $\|P\nabla f(x_k)\|_2 \leq \kappa_{\nabla f}$ needs to be replaced by $\|PG_k\|_2 \leq \kappa_{\nabla f} + \sigma$ for all $k \in \N{}$, which follows under Assumption~\ref{ass.stochastic}.  This is the reason that $\beta^\sigma$ in Parameter Rule~\ref{par.stochastic} involves an additional $\sigma$ term, as opposed to $\beta$ in Lemma~\ref{lem.delta}, which does not.  This affects the denominators in the expression for $\gamma_{k,i,\min}$ for all $(k,i) \in \N{} \times [m]$.  Otherwise, $\gamma_{k,i,\min}$ is consistent between Lemma~\ref{lem.delta} and Parameter Rule~\ref{par.stochastic}.

\blemma\label{lem.delta_stochastic}
  Suppose that Assumptions~\ref{ass.main}, \ref{ass.c}, \ref{ass.bound}, \ref{ass.stochastic}, and \ref{ass.H} hold and that Algorithm~\ref{alg.slip} is run with Parameter Rule~\ref{par.stochastic}.  Then, for all $k \in \N{}$, the choice of $\Gamma_k \in (0,\gamma_{k,\max}]$ in Parameter Rule~\ref{par.stochastic} ensures that $[X_k,X_k + \Gamma_k \alpha_k D_k] \subseteq \Ncal(\theta_k)$, i.e., the condition in Step~\ref{step.gamma} of Algorithm~\ref{alg.slip} holds.  Consequently, under Parameter Rule~\ref{par.stochastic}, Algorithm~\ref{alg.slip} is well defined in the sense that it will generate an infinite sequence of iterates.  Moreover, Parameter Rule~\ref{par.stochastic} ensures that $\|D_k\|_2 \leq \beta^\sigma$ and $\Gamma_k \in [\gamma_{k,\min}, \gamma_{k,\max}]$ for all $k \in \N{}$.
\elemma

Now let us turn to results for the stochastic setting that are distinct from those for the deterministic setting.  Our first main goal is to prove an upper bound on the expected decrease in the shifted barrier-augmented function that occurs with each iteration.  Toward this end, let us begin with the following preliminary bound.

\blemma\label{lem.decrease_stochastic_prelim}
  Suppose that Assumptions~\ref{ass.main}, \ref{ass.c}, \ref{ass.bound}, \ref{ass.stochastic}, and \ref{ass.H} hold and that Algorithm~\ref{alg.slip} is run with Parameter Rule~\ref{par.stochastic}.  Then, for all $k \in \N{}$, one finds that
  \begin{align*}
    &\ \barrier(X_{k+1}, \mu_{k+1}) - \barrier(X_k, \mu_k) \\
    \leq&\ -\Gamma_k \alpha_k \overline\lambda^{-1} \|P\nabla_x \barrier (X_k, \mu_k)\|_2^2 + \thalf \Gamma_k^2 \alpha_k^2 L_k \underline\lambda^{-2} \|PQ_k\|_2^2 \\
    &\ + \Gamma_k \alpha_k \nabla_x \barrier (X_k, \mu_k)^T Z(Z^TH_kZ)^{-1} Z^T (\nabla_x \barrier (X_k, \mu_k) - Q_k).
  \end{align*}
\elemma
\bproof
  Consider arbitrary $k \in \N{}$.  By Lemmas~\ref{lem.lipschitz} and \ref{lem.delta_stochastic}, \eqref{eq.equal_gradients}, \eqref{eq.d_conditions}, \eqref{eq.DZQ}, and the fact that $\|Z^T\nabla_x \barrier(x_k,\mu_k)\|_2 = \|P \nabla_x \barrier(x_k,\mu_k)\|_2$, it follows with $L_k \in \R{}_{>0}$ defined as in Parameter Rule~\ref{par.stochastic} (and \ref{par.deterministic}) that
  \begin{align*}
    &\ \barrier(X_{k+1}, \mu_k) - \barrier(X_k, \mu_k) \\
    \leq&\ \nabla_x \barrier (X_k, \mu_k)^T (X_{k+1} - X_k) + \thalf L_k \|X_{k+1} - X_k\|_2^2 \\
    =&\    \Gamma_k \alpha_k \nabla_x \barrier (X_k, \mu_k)^T D_k + \thalf \Gamma_k^2 \alpha_k^2 L_k \|D_k\|_2^2 \\
    =&\    -\Gamma_k \alpha_k \nabla_x \barrier (X_k, \mu_k)^T Z(Z^TH_kZ)^{-1} Z^T Q_k + \thalf \Gamma_k^2 \alpha_k^2 L_k \|Z(Z^TH_kZ)^{-1} Z^TQ_k\|_2^2 \\
    =&\    -\Gamma_k \alpha_k \|Z^T\nabla_x \barrier (X_k, \mu_k)\|_{(ZH_kZ)^{-1}}^2 + \thalf \Gamma_k^2 \alpha_k^2 L_k \|(Z^TH_kZ)^{-1}Z^TQ_k\|_2^2 \\
    &\ +   \Gamma_k \alpha_k \nabla_x \barrier (X_k, \mu_k)^T Z(Z^TH_kZ)^{-1} Z^T (\nabla_x \barrier (X_k, \mu_k) - Q_k) \\
    \leq&\    -\Gamma_k \alpha_k \overline\lambda^{-1} \|P\nabla_x \barrier (X_k, \mu_k)\|_2^2 + \thalf \Gamma_k^2 \alpha_k^2 L_k \underline\lambda^{-2} \|PQ_k\|_2^2 \\
    &\ +   \Gamma_k \alpha_k \nabla_x \barrier (X_k, \mu_k)^T Z(Z^TH_kZ)^{-1} Z^T (\nabla_x \barrier (X_k, \mu_k) - Q_k).
  \end{align*}
  Finally, since $\mu_{k+1} \leq \mu_k$ implies $\barrier (X_{k+1}, \mu_{k+1}) \leq \barrier (X_{k+1}, \mu_k)$ (see the discussion following~\eqref{eq.barrier_bounded}), the desired conclusion follows. \qed
\eproof

Lemma~\ref{lem.decrease_stochastic_prelim} shows that the change in the shifted barrier-augmented function with each iteration is bounded above by the sum of a negative term and two terms that may be considered noise terms.  The goal of our next two lemmas is to bound these noise terms.  We start with the second term on the right-hand side in Lemma~\ref{lem.decrease_stochastic_prelim}.

\blemma\label{lem.noise_bnd}
  Suppose that Assumptions~\ref{ass.main}, \ref{ass.c}, \ref{ass.bound}, \ref{ass.stochastic}, and \ref{ass.H} hold and that Algorithm~\ref{alg.slip} is run with Parameter Rule~\ref{par.stochastic}.  Then, for all $k \in \N{}$, one finds that
  \bequationNN
    \thalf \Gamma_k^2 \alpha_k^2 L_k \underline\lambda^{-2} \|PQ_k\|_2^2 \leq \tfrac34 \Gamma_k^2 \alpha_k^2 L_k \underline\lambda^{-2} \|P \nabla_x \barrier(X_k, \mu_k)\|_2^2 + \tfrac32 \alpha_k^2 L_k \underline\lambda^{-2} \sigma^2.
  \eequationNN
\elemma
\bproof
  Consider arbitrary $k \in \N{}$. Observe that, for any $(a,b) \in \R{n} \times \R{n}$,
  \bequationNN
    0 \leq \|P(\thalf a - b)\|_2^2 = \tfrac14 \|Pa\|_2^2 - a^TP^TPb + \|Pb\|_2^2,
  \eequationNN
  from which it follows that
  \bequationNN
    \|P(a+b)\|_2^2 = \|Pa\|_2^2 + 2a^TP^TPb + \|Pb\|_2^2 \leq \tfrac{3}{2} \|Pa\|_2^2 + 3 \|Pb\|_2^2.
  \eequationNN
  Consequently, with \eqref{eq.equal_gradients} it holds that
  \begin{align*}
    \thalf \|PQ_k\|_2^2
    =&\ \thalf \|P(\nabla_x \barrier(X_k, \mu_k) + Q_k - \nabla_x \barrier(X_k, \mu_k))\|_2^2 \\ 
    \leq &\ \tfrac34 \|P \nabla_x \barrier(X_k, \mu_k)\|_2^2 + \tfrac32 \|P(Q_k - \nabla_x \barrier(X_k, \mu_k))\|_2^2 \\ 
    \leq &\ \tfrac34 \|P \nabla_x \barrier(X_k, \mu_k)\|_2^2 + \tfrac32 \sigma^2.
  \end{align*}
  Therefore, since $\Gamma_k \leq 1$ under Parameter Rule~\ref{par.stochastic}, the proof is complete. \qed
\eproof

The next lemma provides an upper bound on a conditional expectation of the last term on the right-hand side of the inequality in Lemma~\ref{lem.decrease_stochastic_prelim}.

\blemma\label{lem.inprod_bnd}
  Suppose that Assumptions~\ref{ass.main}, \ref{ass.c}, \ref{ass.bound}, \ref{ass.stochastic}, and \ref{ass.H} hold and that Algorithm~\ref{alg.slip} is run with Parameter Rule~\ref{par.stochastic}.  Then, for all $k \in \N{}$, one finds that
  \bequationNN
    \E [\Gamma_k \alpha_k \nabla_x \barrier (X_k, \mu_k)^T Z(Z^TH_kZ)^{-1} Z^T (\nabla_x \barrier (X_k, \mu_k) - Q_k) | \Gcal_k] \leq \gamma_{k,\buff} \alpha_k B^\sigma,
  \eequationNN
  where
  \bequationNN
    B^\sigma := \underline\lambda^{-1} \(\kappa_{\nabla f} + \frac{\mu_1}{\theta_0} \sum_{j \in [m]} \kappa_{\nabla c_j}\) \sigma.
  \eequationNN
\elemma
\bproof
  Consider arbitrary $k \in \N{}$.  The inner product in the expected value of interest may be nonnegative or nonpositive.  Hence, let us invoke the Law of Total Expectation in order to handle each case separately.  Let $\Ical_k$ be the event that $I_k \geq 0$, where $I_k := \nabla_x \barrier (X_k, \mu_k)^T Z(Z^TH_kZ)^{-1} Z^T (\nabla_x \barrier (X_k, \mu_k) - Q_k)$, and let $\Ical_k^c$ be the complementary event that $I_k < 0$.  By Parameter Rule~\ref{par.stochastic}, the fact that $\E[I_k | \Gcal_k] = 0$ under Assumption~\ref{ass.stochastic}, and Lemma~\ref{lem.delta_stochastic}, it follows that
  \begin{align}
    &\ \E[ \Gamma_k \alpha_k I_k | \Gcal_k] \nonumber \\
    =&\ \E[ \Gamma_k \alpha_k I_k | \Gcal_k \wedge \Ical_k] \P[\Ical_k | \Gcal_k] + \E[ \Gamma_k \alpha_k I_k | \Gcal_k \wedge \Ical_k^c] \P[\Ical_k^c | \Gcal_k] \nonumber \\
    \leq&\ \gamma_{k,\max} \alpha_k \E[I_k | \Gcal_k \wedge \Ical_k] \P[\Ical_k | \Gcal_k] + \gamma_{k,\min} \alpha_k \E[I_k | \Gcal_k \wedge \Ical_k^c] \P[\Ical_k^c | \Gcal_k] \nonumber \\
    \leq&\ \gamma_{k,\min} \alpha_k (\E[I_k | \Gcal_k \wedge \Ical_k] \P[\Ical_k | \Gcal_k] + \E[I_k | \Gcal_k \wedge \Ical_k^c] \P[\Ical_k^c | \Gcal_k]) \nonumber \\
    &\ + \gamma_{k,\buff} \alpha_k \E[I_k | \Gcal_k \wedge \Ical_k] \P[\Ical_k | \Gcal_k] \nonumber \\
    =&\ \gamma_{k,\buff} \alpha_k \E[I_k | \Gcal_k \wedge \Ical_k] \P[\Ical_k | \Gcal_k]. \label{eq.noisynoise}
  \end{align}
  Furthermore, by Assumptions~\ref{ass.main}--\ref{ass.H}, the Cauchy--Schwarz inequality, the triangle inequality, $\|P\|_2 \leq 1$, and the fact that $\|Z^Tq\|_2 = \|Pq\|_2$ for any $q \in \R{n}$ one has
  \begin{align*}
    &\ \E[I_k | \Gcal_k \wedge \Ical_k] \P[\Ical_k | \Gcal_k] \\ 
    =&\ \E[ \nabla_x \barrier (X_k, \mu_k)^T Z(Z^TH_kZ)^{-1} Z^T (\nabla_x \barrier (X_k, \mu_k) - Q_k) | \Gcal_k \wedge \Ical_k] \P[\Ical_k | \Gcal_k] \\
    \leq&\ \underline\lambda^{-1} \E[ \|P \nabla_x \barrier (X_k, \mu_k)\|_2 \|P (\nabla_x \barrier (X_k, \mu_k) - Q_k)\|_2 | \Gcal_k \wedge \Ical_k] \P[\Ical_k | \Gcal_k] \\
    \leq &\ \underline\lambda^{-1} \(\kappa_{\nabla f} + \frac{\mu_1}{\theta_0} \sum_{j \in [m]} \kappa_{\nabla c_j}\) \sigma.
  \end{align*}
  Combining this inequality with \eqref{eq.noisynoise} yields the desired conclusion.
  \qed
\eproof

Combining the preliminary bound in Lemma~\ref{lem.decrease_stochastic_prelim} with the previous two lemmas, we are now prepared to prove an upper bound on the expected decrease in the shifted barrier-augmented function that occurs with each iteration.

\blemma\label{lem.decrease_stochastic}
  Suppose that Assumptions~\ref{ass.main}, \ref{ass.c}, \ref{ass.bound}, \ref{ass.stochastic}, and \ref{ass.H} hold and that Algorithm~\ref{alg.slip} is run with Parameter Rule~\ref{par.stochastic}.  Then, for all $k \in \N{}$, one finds that
  \begin{align*}
    &\ \E[\barrier (X_{k+1},\mu_{k+1}) | \Gcal_k] - \barrier(X_k, \mu_k) \\
    \leq&\ - \tfrac14 \gamma_{k,\min} \alpha_k \overline\lambda^{-1} \|P\nabla_x \barrier (X_k, \mu_k)\|_2^2 + \tfrac32 \alpha_k^2 L_k \underline\lambda^{-2} \sigma^2 + \gamma_{k,\buff} \alpha_k B^\sigma,
  \end{align*}
  where $B^\sigma$ is defined as in Lemma~\ref{lem.inprod_bnd}.
\elemma
\bproof
  Consider arbitrary $k \in \N{}$. One finds from Lemmas \ref{lem.decrease_stochastic_prelim}, \ref{lem.noise_bnd}, and \ref{lem.inprod_bnd} that
  \begin{align*}
    &\ \E[\barrier (X_{k+1},\mu_{k+1}) | \Gcal_k] - \barrier(X_k, \mu_k) \\ 
    \leq &\ -\Gamma_k \alpha_k (\overline\lambda - \tfrac34 \Gamma_k \alpha_k L_k \underline\lambda^{-2}) \|P\nabla_x \barrier (X_k, \mu_k)\|_2^2 + \tfrac32 \alpha_k^2 L_k \underline\lambda^{-2} \sigma^2 + \gamma_{k,\buff} \alpha_k B^\sigma.
  \end{align*}
  Then, under Parameter Rule~\ref{par.stochastic}, one has $\Gamma_k \leq 1$ and $k^{t_\alpha} \leq 1$, so
  \bequationNN
    \alpha_k = k^{t_\alpha} \underline\lambda^2 \overline\lambda^{-1} L_k^{-1} \implies \overline\lambda - \tfrac34 \Gamma_k \alpha_k L_k \underline\lambda^{-2} = \overline\lambda^{-1} - \tfrac34 \overline\lambda^{-1} \Gamma_k k^{t_\alpha} \geq \tfrac14 \overline\lambda^{-1},
  \eequationNN
  which when combined with the prior conclusion completes the proof.
  \qed
\eproof

Using this sufficient decrease result and looking further into the selections in Parameter Rule~\ref{par.stochastic}, we obtain the following result.

\blemma\label{lem.stochastic_combined}
  Suppose that Assumptions~\ref{ass.main}, \ref{ass.c}, \ref{ass.bound}, \ref{ass.stochastic}, and \ref{ass.H} hold and that Algorithm~\ref{alg.slip} is run with Parameter Rule~\ref{par.stochastic}.  Then, $(\nu,\varsigma) \in \R{}_{>0} \times \R{}_{>0}$ exists such that, for all $k \in \N{}$,
  \begin{multline*}
    \E[ \barrier(X_{k+1},\mu_{k+1}) | \Gcal_k] - \barrier(X_k,\mu_k) \\ \leq - \nu k^{t + t_\alpha} \|P \nabla_x \barrier(X_k,\mu_k) \|_2^2 + \varsigma k^{\max\{2t + t_\alpha,t + 2t_\alpha\}}.
  \end{multline*}
\elemma
\bproof
  We prove the result by building on Lemma~\ref{lem.decrease_stochastic}.  By Parameter Rule~\ref{par.stochastic}, there exist $\kappa \in \R{}_{>0}$, $\underline{L} \in \R{}_{>0}$, and $L \in \R{}_{>0}$ such that for all $k \in \N{}$ one finds
  \bequationNN
    \gamma_{k,\min} \alpha_k = \min \left\{ \alpha_k, \frac{-\kappa + \sqrt{\kappa^2 + 2 \underline{L} (\eta \mu_k - \theta_k)}}{\beta^\sigma L}, \frac{\underline\eta - \theta_k}{\beta^\sigma L} \right\}.
  \eequationNN
  Since this is the same type of lower bound as in \eqref{def:betak-alphak}, the same argument in the proof of Theorem~\ref{th.deterministic} shows that there exists $\nu \in \R{}_{>0}$ such that $\tfrac14 \gamma_{k,\min} \alpha_k \overline\lambda^{-1} \geq \nu k^{t + t_\alpha}$ for all $k \in \N{}$.  Hence, by Lemma~\ref{lem.decrease_stochastic}, all that remains is to prove that there exists $\varsigma \in \R{}_{>0}$ such that for all $k \in \N{}$ one finds that
  \bequation\label{eq.zeb}
    \tfrac32 \alpha_k^2 L_k \underline\lambda^{-2} \sigma^2 + \gamma_{k,\buff} \alpha_k B^\sigma \leq \varsigma k^{\max\{2t + t_\alpha,t + 2t_\alpha\}}.
  \eequation
  Starting with the first term in \eqref{eq.zeb}, one finds that
  \bequationNN
    \tfrac32 \alpha_k^2 L_k \underline\lambda^{-2} \sigma^2 = \tfrac32 \alpha_k (\alpha_k L_k \underline\lambda^{-2}) \sigma^2 = \tfrac32 \underline\lambda^2 \overline\lambda^{-2} \sigma^2 \frac{k^{2t_\alpha}}{L_k},
  \eequationNN
  where with $\upsilon := \sum_{i=1}^m \(L_{c_i} \kappa_{\nabla c_i} + \kappa_{c_i} L_{\nabla c_i}\)$ one finds that
  \bequationNN
    \frac{k^{2t_\alpha}}{L_k} = \frac{k^{2t_\alpha}}{L_{\nabla f} + \mu_1 k^t \theta_0^{-2} k^{-t} (k+1)^{-t} \upsilon} \leq \frac{k^{2t_\alpha}}{L_{\nabla f} + \mu_1 \theta_0^{-2} \upsilon k^{-t}} \leq \frac{k^{t + 2t_\alpha}}{\mu_1 \theta_0^{-2} \upsilon}.
  \eequationNN
  As for the second term in \eqref{eq.zeb}, one similarly finds that
  \bequationNN
    \gamma_{k,\buff} \alpha_k B^\sigma \leq \gamma_{\buff} k^t k^{t_\alpha} \(\frac{\underline\lambda^2 B^\sigma}{\overline\lambda^2 L_k}\) \leq \(\frac{\gamma_{\buff} \underline\lambda^2 B^\sigma}{\mu_1 \theta_0^{-2} \upsilon \overline\lambda^2}\) k^{2t + t_\alpha}.
  \eequationNN
  Combining these bounds yields \eqref{eq.zeb} for some $\varsigma \in \R{}_{>0}$.
  \qed
\eproof

Our final theorem, presented next, now follows similarly to Theorem~3.16 in~\cite{CKRW23}.  We provide a complete proof for the sake of completeness.

\btheorem\label{th.main}
  Suppose that Assumptions~\ref{ass.main}, \ref{ass.c}, \ref{ass.bound}, \ref{ass.stochastic}, and \ref{ass.H} hold and that Algorithm~\ref{alg.slip} is run with Parameter Rule~\ref{par.stochastic}.  Then,
  \bequationNN
    \liminf_{k \to \infty} \|P \nabla_x \phi(X_k,\mu_k)\|_2^2 = 0\ \ \text{almost surely}.
  \eequationNN
  Consequently, if in a given run there exists $\Kcal \subseteq \N{}$ with $|\Kcal| = \infty$ such that
  \bitemize
    \renewcommand\labelitemi{$\bullet$}
    \item $\{P\nabla_x \phi(x_k,\mu_k)\}_{k\in\Kcal} \to 0$,
    \item $\{x_k\}_{k\in\Kcal} \to \xbar$ for some $\xbar \in \Fcal$, and
    \item at $\xbar$ the linear independence constraint qualification (LICQ) holds with respect to problem~\eqref{prob.opt} in the sense that with $\bar\Acal := \{i \in [m] : c_i(\xbar) = 0\}$ the columns of $A^T$ combined with the vectors in $\{\nabla c_i(\xbar)\}_{i\in\bar\Acal}$ form a linearly independent set,
  \eitemize
  then $\xbar$ is a KKT point for problem~\eqref{prob.opt} in the sense that there exists a pair of Lagrange multipliers $(\ybar,\zbar) \in \R{l} \times \R{m}$ such that $(\xbar,\ybar,\zbar)$ satisfies \eqref{eq.KKT}.
\etheorem
\bproof
  It follows from Lemma~\ref{lem.stochastic_combined} and the Law of Total Expectation that there exists $(\nu,\varsigma) \in \R{}_{>0} \times \R{}_{>0}$ such that, for all $k \in \N{}$, one finds that
  \begin{multline*}
    \E[ \barrier(X_{k+1},\mu_{k+1})] - \E[\barrier(X_k,\mu_k)] \\ \leq - \nu k^{t + t_\alpha} \E[\|P \nabla_x \barrier(X_k,\mu_k) \|_2^2] + \varsigma k^{\max\{2t + t_\alpha,t + 2t_\alpha\}}.
  \end{multline*}
  Consider arbitrary $K \in \N{}$.  Summing over $k \in [K]$ yields
  \begin{align*}
    f_{\inf} - \E[\barrier(x_1,\mu_1)]
      &\leq \E[\barrier(X_{K+1},\mu_{K+1})] - \E[\barrier(x_1,\mu_1)] \\
      &\leq - \nu \sum_{k=1}^{K} k^{t + t_\alpha} \E[\|P \nabla_x \barrier(X_k,\mu_k)\|_2^2] + \varsigma \sum_{k=1}^{K} k^{\max\{2t + t_\alpha,t+2t_\alpha\}}.
  \end{align*}
  After rearrangement, this yields
  \bequationNN
    \sum_{k=1}^{K} k^{t + t_\alpha} \E[\|P\nabla_x \barrier(X_k,\mu_k)\|_2^2] \leq \tfrac{1}{\nu} (\E[\barrier(x_1,\mu_1)] - f_{\inf}) + \tfrac{\varsigma}{\nu} \sum_{k=1}^{K} k^{\max\{2t + t_\alpha ,t + 2t_\alpha\}}.
  \eequationNN
  Since $2t + t_\alpha \in (-\infty,-1)$ and $t + 2t_\alpha \in (-\infty,-1)$, the right-hand side of this inequality converges to a finite limit as $K \to \infty$.  Furthermore, by $\sum_{k=1}^\infty k^{t + t_\alpha} = \infty$, nonnegativity of $\|P \nabla_x \barrier (X_k, \mu_k)\|_2^2$, and Fatou's lemma, one almost surely has
  \bequationNN
    0 = \liminf_{k \to \infty} \E[\|P \nabla_x \barrier (X_k, \mu_k)\|_2^2] \geq \E\left[\liminf_{k \to \infty} \|P \nabla_x \barrier (X_k, \mu_k)\|_2^2 \right] = 0.
  \eequationNN
  Now consider $\Phi := \liminf_{k \to \infty} \|P \nabla_x \barrier (X_k, \mu_k)\|_2^2$, the expectation of which has been shown above to be zero.  Nonnegativity of $\|P \nabla_x \barrier (X_k,\mu_k)\|_2^2$ for all $k \in \N{}$ and the Law of Total Expectation imply that $0 = \E[\Phi] \geq \P[\Phi > 0] \E[\Phi | \Phi > 0]$, so $0 = \P[\Phi > 0] = \P[\liminf_{k \to \infty} \|P \nabla_x \barrier (X_k, \mu_k)\|_2^2 > 0]$.  This is the first desired conclusion of the theorem.  The second desired conclusion, on the other hand, follows by using the same argument as in the proof of Theorem~\ref{th.deterministic}.
  \qed
\eproof

%*********
% Section
%*********
\section{Discussion of Assumptions~\ref{ass.c} and \ref{ass.H}}\label{sec.assumption}

Our aim in this section is to discuss Assumptions~\ref{ass.c} and \ref{ass.H}.  In particular, first for the deterministic setting, we provide (a) example problems for which Assumption~\ref{ass.c} can be shown to hold or cannot be shown to hold, (b) further explanation about why Assumption~\ref{ass.c} can be viewed as a combination of a nondegeneracy assumption and an assumption about $\mu_1/\theta_0$ being sufficiently large, and (c) a procedure for computing~$d_k$ to satisfy~\eqref{eq.d_conditions} when such a direction exists.  We follow this with a discussion about the combination of Assumptions~\ref{ass.c} and \ref{ass.H} for the stochastic setting.  We emphasize upfront that, generally speaking, \emph{a straightforward procedure for computing $d_k$ to satisfy the conditions of our convergence analysis $($for either the deterministic or stochastic settings$)$ may not be available in practice.}  For example, for one thing, determining a sufficiently large value for the ratio $\mu_1/\theta_0$ may rely on problem-specific constants that are not known in advance of a run of the algorithm.  That being said, in this section we provide guidance for how to compute $d_k$ in practice such that, along with parameter tuning (as is common for stochastic-gradient-based methods in practice), one obtains good practical performance.

Throughout this section, for the sake of clarity, we drop the subscript notation that indicates the iteration index of an algorithm---rather, in this section, a subscript only refers to an index of a component of a vector.  In particular, in this section, we denote $x \equiv x_k$ (which in turn means, e.g., that $x_1$ refers to the first component of the vector $x$), $\mu \equiv \mu_k$, $\theta \equiv \theta_{k-1}$, $q \equiv q_k$, and $d \equiv d_k$.

Let us begin with some example problems for the deterministic setting.  Firstly, recall that under Assumption~\ref{ass.main} a strictly feasible point must exist.  This means that one can rule out problems with, e.g., the affine constraint $x_1 = 0$ and the inequality constraint $x_1 \leq 0$, or any problems with such a combination of constraints that implies the lack of a strictly feasible point.  Secondly, observe that the existence or not of a direction satisfying \eqref{eq.d_conditions} depends on the feasible region reduced to $\Null(A)$.  This means that one can distinguish cases of whether \eqref{eq.d_conditions} holds or not by presuming that one has already restricted attention to this reduced space and consider only the geometry of the feasible region with respect to inequality constraints where the ambient space of the variables is this reduced space.

Let us now present an example for which Assumption~\ref{ass.c} can be shown to hold.  In particular, the following example shows that there exists a value of the tuple $(\mu,\theta,\eta,\underline\eta,\underline\zeta,\overline\zeta,\zeta)$ under the ranges specified by Algorithm~\ref{alg.slip} such that at all $x \in \Fcal_{<0}$ there exists $d$ satisfying \eqref{eq.d_conditions}.  The main aspect of the example is that $\mu$ must be sufficiently large relative to $\theta$.  We remark that the case $a \in \R{}_{>0}$ that is considered in the example is the more challenging case to consider.  Indeed, if $a \in \R{}_{\leq0}$, then, using similar arguments, the parameter requirements are less restrictive than in the following example.  We expound on this claim further after the example.

\begin{example}\label{ex.2.dimensional}
  To start, let $(\mu,\theta,\eta,\underline\eta,\underline\zeta,\overline\zeta,\zeta)$ be arbitrary positive constants under the ranges specified by Algorithm~\ref{alg.slip}.  (A further restriction on the relationship between $\mu$ and $\theta$ is developed during the discussion of this example, when it is needed.)  Consider a two-dimensional problem where $f(x) = v^Tx$ for some $v \in \R{2}$, $c_1(x) = -x_1$, $c_2(x) = ax_1 - x_2$ for some $a \in \R{}_{>0}$, and $g = \nabla f(x) = v$ (see~\eqref{eq.gradient}), so
  \bequation\label{ex.1.q}
    q = v + \mu \(\frac{1}{-c_1(x)}\) \nabla c_1(x) + \mu \(\frac{1}{-c_2(x)}\) \nabla c_2(x).
  \eequation
  Defining $\Acal(x) := \{i \in \{1,2\} : -\eta \mu < c_i(x) \leq -\theta\}$, the conditions~\eqref{eq.d_conditions} reduce to
  \bsubequations\label{ex.1}
  \begin{align}
    \underline\zeta \|q\|_2 \leq \|d\|_2 &\leq \overline\zeta \|q\|_2, \label{ex.1.1} \\
    -q^Td &\geq \zeta \|q\|_2 \|d\|_2, \label{ex.1.2} \\ \text{and}\ \ 
    \nabla c_i(x)^Td &\leq -\thalf \underline\eta \|d\|_2\ \ \text{for all}\ \ i \in \Acal(x). \label{ex.1.3}
  \end{align}
  \esubequations
  Since the gradient vectors $\nabla f(x)$, $\nabla c_1(x)$, and $\nabla c_2(x)$ are constant over $x \in \Fcal_{<0}$ and the latter two gradient vectors are nonzero, it follows from~\eqref{ex.1.q} that for any $x \in \Fcal_{<0}$ one finds that $q/\|q\|_2 \to q_{\lim}/\|q_{\lim}\|_2$ as $\mu \to \infty$, where $q_{\lim} = \nabla c_1(x)/(-c_1(x)) + \nabla c_2(x)/(-c_2(x))$.  At the same time, for any pair of constants $(\chi_1,\chi_2) \in \R{}_{>0} \times \R{}_{>0}$ and with the search direction defined by
  \bequation\label{ex.1.d}
    d = -v - \mu \(\frac{\chi_1}{-c_1(x)}\) \nabla c_1(x) - \mu \(\frac{\chi_2}{-c_2(x)}\) \nabla c_2(x),
  \eequation
  one has for any $x \in \Fcal_{<0}$ that $d/\|d\|_2 \to d_{\lim}/\|d_{\lim}\|_2$ as $\mu \to \infty$, where the limiting direction is $d_{\lim} = (\chi_1/c_1(x)) \nabla c_1(x) + (\chi_2/c_2(x)) \nabla c_2(x)$.
  
  Now suppose that the constants $(\mu,\theta,\eta)$ and point $x \in \Fcal_{<0}$ are such that both constraints are nearly active at $x$, i.e., $\Acal(x) = \{1,2\}$.  Since one has $\chi_3 := (-\nabla c_1(x)^T\nabla c_2(x))/(\|\nabla c_1(x)\|_2 \|\nabla c_2(x)\|_2) = a/\sqrt{1 + a^2} \in (0,1)$, it follows that
  \begin{align*}
    \nabla c_1(x)^Td_{\lim} &= \chi_1 \underbrace{\frac{\|\nabla c_1(x)\|_2^2}{c_1(x)}}_{< 0} + \chi_2 \chi_3 \underbrace{\frac{\|\nabla c_1(x)\|_2 \|\nabla c_2(x)\|_2}{-c_2(x)}}_{> 0} \\ \text{and}\ \ 
    \nabla c_2(x)^Td_{\lim} &= \chi_1 \chi_3 \underbrace{\frac{\|\nabla c_1(x)\|_2 \|\nabla c_2(x)\|_2}{-c_1(x)}}_{> 0} + \chi_2 \underbrace{\frac{\|\nabla c_2(x)\|_2^2}{c_2(x)}}_{<0}.
  \end{align*}
  Now suppose $\chi_1 = c_1(x)\|\nabla c_2(x)\|_2/(c_2(x)\|\nabla c_1(x)\|_2)$ and $\chi_2 = 1$.  (The effect of these choices is to normalize the contributions of the two terms in the definition of $d_{\lim}$.)  These choices and the fact that $\eta = \psi \theta / \mu$ for some $\psi \in \R{}_{>1}$ implies that
  \begin{align*}
    \nabla c_1(x)^Td_{\lim} &< - \( \frac{1 - \chi_3}{\psi \theta} \) \|\nabla c_1(x)\|_2 \|\nabla c_2(x)\|_2 \\ \text{and}\ \ 
    \nabla c_2(x)^Td_{\lim} &< - \( \frac{1 - \chi_3}{\psi \theta} \) \|\nabla c_2(x)\|_2^2.
  \end{align*}
  On the other hand, one has that $\|d_{\lim}\|_2 \leq 2 \|\nabla c_2(x)\|_2 / \theta$.  Thus, one may conclude that there exists $\chi_4 \in \R{}_{>0}$ such that $\nabla c_1(x)^Td_{\lim} \leq - \chi_4 \|d_{\lim}\|_2$ and $\nabla c_2(x)^Td_{\lim} \leq - \chi_4 \|d_{\lim}\|_2$ for all $x \in \Fcal_{<0}$ such that both constraints are nearly active.
  
  On the other hand, with $d_{\lim}$ defined as in the previous paragraph (i.e., with $\chi_1 = c_1(x) \|\nabla c_2(x)\|_2/(c_2(x)\|\nabla c_1(x)\|_2)$ and $\chi_2 = 1$), it follows that
  \begin{align*}
    &\ -q_{\lim}^Td_{\lim} \\
    =& \(\frac{1}{c_1(x)} \nabla c_1(x) + \frac{1}{c_2(x)} \nabla c_2(x)\)^T \(\frac{\chi_1}{c_1(x)} \nabla c_1(x) + \frac{\chi_2}{c_2(x)} \nabla c_2(x)\) \\
    =&\ \frac{1}{c_1(x)c_2(x)} \|\nabla c_1(x)\|_2 \|\nabla c_2(x)\|_2 + \frac{1}{c_2(x)^2} \|\nabla c_2(x)\|_2^2 \\
    &\ + \( \frac{1}{c_1(x) c_2(x)} + \frac{\|\nabla c_2(x)\|_2}{c_2(x)^2 \|\nabla c_1(x)\|_2} \) \nabla c_1(x)^T\nabla c_2(x) \\
    =&\ \frac{1 - \chi_3}{c_1(x)c_2(x)} \|\nabla c_1(x)\|_2 \|\nabla c_2(x)\|_2 + \frac{1 - \chi_3}{c_2(x)^2} \|\nabla c_2(x)\|_2^2 \\
    \geq&\ \frac{1 - \chi_3}{\theta^2} (\|\nabla c_1(x)\|_2 \|\nabla c_2(x)\|_2 + \|\nabla c_2(x)\|_2^2).
  \end{align*}
  On the other hand, $\|d_{\lim}\|_2 \leq 2 \|\nabla c_2(x)\|_2 / \theta$ and $\|q_{\lim}\|_2 \leq 2 \|\nabla c_2(x)\|_2 / \theta$.  Thus, one may conclude that there exists $\chi_5 \in \R{}_{>0}$ such that $-q_{\lim}^Td_{\lim} \geq \chi_5 \|q_{\lim}\|_2 \|d_{\lim}\|_2$ for all $x \in \Fcal_{<0}$ such that both constraints are nearly active.
  
  The arguments in the preceding two paragraphs involve the limiting directions $d_{\lim}$ and $q_{\lim}$, where it was shown that there exist positive constants $\chi_4$ and $\chi_5$---that are uniform with respect to $x \in \Fcal_{<0}$---such that conditions of the form \eqref{ex.1.2}--\eqref{ex.1.3} hold for these limit directions.  All that remains is to observe that, since these are the limiting directions for any given $x \in \Fcal_{<0}$ as $\mu \to \infty$, such conditions and \eqref{ex.1.1} hold for $(d,q)$ for some positive constants as long as $\mu$ is sufficiently large relative to $\theta$.  All together from our observations, it follows that with $\mu$ sufficiently large relative to $\theta$, $\eta \in (\theta/\mu,1)$, $x \in \Fcal_{<0}$ such that both constraints are nearly active, $\underline\eta$ sufficiently small, $\underline\zeta$ sufficiently small, $\overline\zeta$ sufficiently large, and $\zeta$ sufficiently small, one finds $d$ in \eqref{ex.1.d} with the aforementioned choices of $\chi_1$ and $\chi_2$ satisfies~\eqref{eq.d_conditions}.
  
  Finally, observe that for any $x \in \Fcal_{<0}$ such that both constraints are not nearly active, the situation is simpler than above.  In such a setting, one does not need condition~\eqref{ex.1.3} to hold for both constraints; it only needs to hold for whichever constraint, if any, is nearly active.  This offers much more flexibility in the choice of~$d$.  Overall, since the constraint gradients are constant over $x \in \Fcal_{<0}$, one can again derive the existence of uniform parameter choices, including that $\mu$ is sufficiently large relative to $\theta$, to ensure that $d$ satisfying \eqref{eq.d_conditions} exists at all such $x$.
  \qed
\end{example}

The analysis in Example~\ref{ex.2.dimensional} can be extended to any problem with polyhedral constraints (for arbitrary $n \in \N{}$ and $m \in \N{}$) as long as at any strictly feasible point in $\Fcal_{<0}$ the interior of the polar cone of the nearly active constraint gradients is nonempty.  Intuitively, this can be understood as follows.  With $\mu$ sufficiently large relative to $\theta$, the gradient of the barrier term corresponding to nearly active constraints \emph{pushes} the vector $-q$ (through more emphasis on the barrier term) to point with the interior of the polar cone of the nearly active constraint gradients.  This means that, as long as a direction $d$ is chosen appropriately to compensate for the different magnitudes of the norms of the constraint gradients, the desired conditions in \eqref{eq.d_conditions} hold.  Admittedly, \emph{the parameter choices that ensure that such a direction $d$ always exists are not straightforward to determine in advance of a run of the algorithm}.  This means that, in practice, the algorithm parameters require tuning and/or adaptive choices to be made, as we discuss further shortly.  To illustrate the situation that we have described in Example~\ref{ex.2.dimensional}, please see Figure~\ref{fig.example1}, on the left.

\bfigure[ht]
  \centering
  \begin{tikzpicture}
  \coordinate (o)  at ( 0.00, 0.00);
  \coordinate [label=left:{$\nabla c_1(x)$}] (g1) at (-1.00, 0.00);
  \coordinate [label=below right:{$\nabla c_2(x)$}] (g2) at ( 1.00,-1.00);
  \coordinate (c1) at ( 0.00, 2.00);
  \coordinate (c2) at ( 2.00, 2.00);
  \coordinate [label=below left:$x$] (x) at ( 0.35, 0.60);
  \coordinate [label=above:{$x-q$}] (xq) at ( 0.55, 2.10);
  \coordinate [label=above:{$x+d$}] (xd) at ( 1.20, 1.80);
  \draw[very thick,dashed,-latex] (o) -- (g1);
  \draw[very thick,dashed,-latex] (o) -- (g2);
  \draw[very thick,-latex] (o) -- (c1);
  \draw[very thick,-latex] (o) -- (c2);
  \draw[very thick,gray,-latex] (x) -- (xq);
  \draw[very thick,gray,-latex] (x) -- (xd);
  \filldraw (x) circle (1.00pt);

  \begin{scope}[shift={(5.0,0.0)}]
  \coordinate (o)  at ( 0.00, 0.00);
  \coordinate [label=left:{$\nabla c_1(x)$}] (g1) at (-1.00, 0.00);
  \coordinate [label=below right:{$\nabla c_2(0)$}] (g2) at ( 1.00, 0.00);
  \coordinate (c1) at ( 0.00, 2.00);
  \coordinate [label=below left:$x$] (x) at ( 0.35, 0.80);
  \coordinate [label=above:{$x-q$}] (xq) at ( 0.55, 2.10);
  \coordinate [label=above:{$x+d$}] (xd) at ( 1.20, 1.80);
  \draw[very thick,dashed,-latex] (o) -- (g1);
  \draw[very thick,dashed,-latex] (o) -- (g2);
  \draw[very thick,-latex] (o) -- (c1);
  \draw[very thick,domain= 0.00: 2.00,-latex] plot({\x},{sqrt(\x)});
  \draw[very thick,gray,-latex] (x) -- (xq);
  \draw[very thick,gray,-latex] (x) -- (xd);
  \filldraw (x) circle (1.00pt);
  \end{scope}
\end{tikzpicture}
  \caption{On the left, an illustration of Example~\ref{ex.2.dimensional} for which Assumption~\ref{ass.c} holds as long as the barrier parameter $\mu$ is sufficiently large relative to $\theta$ (recall that $\mu_k/\theta_{k-1}$ is constant in the algorithm) amongst other parameter choices.  Since $\mu$ is sufficiently large relative to $\theta$, it follows that a direction $d$ pointing into the interior of the polar cone of nearly active constraint gradients is also one of sufficient decrease for the barrier-augmented objective function.  On the right, an illustration of Example~\ref{ex.2.dim.failure} for which Assumption~\ref{ass.c} fails to hold since, as $x$ approaches the origin from within the feasible region, there is no minimum value for the ratio $\mu/\theta$ such that a direction into the interior of the polar cone of nearly active constraint gradients satisfies~\eqref{eq.d_conditions}.}
  \label{fig.example1}
\efigure

Now let us briefly describe an example for which we are unable to show that Assumption~\ref{ass.c} holds with uniform parameter choices; see Figure~\ref{fig.example1}, on the right.

\begin{example}\label{ex.2.dim.failure}
  Consider a two-dimensional problem where $f(x) = v^Tx$ for some $v \in \R{2}$, $c_1(x) = -x_1$, $c_2(x) = x_1 - x_2^2$, and $g = \nabla f(x) = v$ (see~\eqref{eq.gradient}).  Defining, as before, $\Acal(x) := \{i \in \{1,2\} : -\eta \mu < c_i(x) \leq -\theta\}$, \eqref{eq.d_conditions} reduces to \eqref{ex.1}.  The situation is similar to that in Example~\ref{ex.2.dimensional}, except that $\nabla c_2(x) = [1\ -2x_2]^T$ is not constant over~$x$.  Attempting to follow a similar analysis as in Example~\ref{ex.2.dimensional}, one finds that the arguments break down since $\chi_3 = -\nabla c_1(x)^T\nabla c_2(x)/(\|\nabla c_1(x)\|_2 \|\nabla c_2(x)\|_2) \to 1$ as $x \to 0$ from within the interior of the feasible region.
  \qed
\end{example}

Intuitively, Example~\ref{ex.2.dim.failure} has the property that the constraint gradients point in opposite directions as $x \to 0$ from within the interior of the feasible region.  This means that, if the algorithm were to compute $x \to 0$ (say, since $0$ is the optimal solution), the barrier terms for the two constraints would push against each other, in the limit.  Consequently, we are not able to show that there exists a uniform ratio $\mu/\theta$ such that a direction $d$ satisfying \eqref{eq.d_conditions} always exists over the feasible region.

One sees in Figure~\ref{fig.example1} our claim that Assumption~\ref{ass.c} is a type of nondegeneracy assumption.  When, at any strictly feasible point, the polar cone of the nearly active constraint gradients is sufficiently \emph{wide}, the situation is favorable, like on the left in Figure~\ref{fig.example1}.  On the other hand, if the polar cone collapses as $x \to \xbar$ for some $\xbar$, like in the situation on the right in Figure~\ref{fig.example1}, then the situation is not favorable since the barrier terms for different constraints push against each other.

We close this section by observing that a procedure for computing a direction to satisfy \eqref{eq.d_conditions}, if one exists, has been suggested in Example~\ref{ex.2.dimensional}.  In particular, at any iterate $x$ one can compute $d$ by (a) determining the set of nearly active constraints, (b) computing weights for the nearly active constraints, as in \eqref{ex.1.d}, such that the contribution from each nearly active constraint gradient is normalized in some sense in order to satisfy \eqref{eq.d_conditions}.  Concretely, one can consider the feasibility problem
\bequationNN
  \text{find}\ \ \{\chi_i\}_{i \in [m]} \in \R{m}_{>0}\ \ \text{such that}\ \ d = g - \mu \sum_{i\in[m]} \frac{\chi_i}{c_i(x)} \nabla c_i(x)\ \ \text{satisfies}\ \ \eqref{eq.d_conditions}.
\eequationNN
However, in practice, we suspect such a procedure not to be worth the required computational effort.  Instead, we suggest that, in practice, a reasonable choice is to initialize $\mu$ and $\theta$ to positive values, and in each iteration compute $d$ by solving a linear system of the form \eqref{eq.linsys}, as we specify in Assumption~\ref{ass.H} for the stochastic setting.  For example, $H_k$ can simply be an identity matrix, or it might be chosen as an identity matrix plus (an approximation of) the Hessian of the barrier function, which may involve second-order derivatives of $\{c_i\}_{i\in[m]}$, but does not require derivatives (of any order) of the objective function $f$.  In any case, if $d$ fails to satisfy \eqref{eq.d_conditions} (specifically, \eqref{eq.d_condition_5}) for some choices of the algorithm parameters, then the algorithm might ``reset'' the barrier parameter to a larger value.  This may be allowed to occur iteratively, at least up to some upper limit for practical purposes.  In the worst case, the algorithm may need to compute smaller values of $\gamma$ than as stated in our theoretical results.  However, at least this is a practical approach that worked well in our experiments.

%*********
% Section
%*********
\section{Numerical Experiments}\label{sec.numerical}

As previously stated, we do not claim that the deterministic version of our algorithm would be competitive with state-of-the-art derivative-based interior-point methods, especially not with such methods that allow infeasible iterates and/or employ second-order derivatives.  On the other hand, for the stochastic setting, we know of no other stochastic interior-point method like ours that can solve generally constrained, smooth optimization problems, meaning that there is no such method in the literature against which our algorithm can be compared on a level playing field.  All of this being said, it is illustrative to demonstrate the practical performance of our method through numerical experiments, both on a well-established and diverse set of test problems and on a couple of relatively straightforward test problems that represent the settings for which the algorithm has been designed.  In this section, we present the results of such numerical experiments.  For our experiments, each run of the algorithm was conducted on a compute node with an 16 AMD Opteron Processor 6128 with 32GB of memory.

First, we conducted an experiment with a large test set of problems in order to demonstrate the broad applicability of our algorithm.  For such a demonstration, in these experiments we considered only the deterministic version of our algorithm.  (Our second and third sets of experiments consider the stochastic version.)  We generated our test set from the CUTEst collection \cite{GoulOrbaToin15}, specifically using the PyCUTEst interface \cite{PyCUTEst2022}.  This collection includes 1315 unconstrained or constrained smooth optimization problems.  To generate the subset of problems that we employed for our experiments, we proceeded as follows.  Firstly, we removed all problems with nonlinear equality constraints, with no inequality constraints, or that were too large to be loaded into memory.  This resulted in a set of 396 problems of the form~\eqref{prob.opt}.  Second, since our algorithm requires a strictly feasible initial point, we ran a so-called Phase~I algorithm (implemented in Python) starting from the initial point provided by PyCUTEst; see Algorithm~\ref{alg.phaseI} in Appendix~\ref{sec.appendix}.  Of the 396 problems on which the Phase~I algorithm was run, 10 encountered function evaluation errors through the interface and 59 resulted in a failure to find a strictly feasible point; these were all removed.  Our test set is the remainder of these problems---327 in total---for which a strictly feasible initial point could be found.  These points were stored as the initial points for this first set of our experiments.

We implemented Algorithm~\ref{alg.slip} in Python.  We refer to the software as \texttt{SLIP}.  For each test problem, for which $x_1$ was determined by Algorithm~\ref{alg.phaseI}, the input parameters for the algorithm were set as follows.  The initial neighborhood parameter was set as $\theta_0 \gets -0.9 \max\{c(x_1)\} > 0$ and the initial barrier parameter was set as $\mu_1 \gets \max\{10^{-1}, 2\theta_0\}$.  The corresponding sequences of neighborhood and barrier parameters were then prescribed such that $\theta_{k-1} \gets \theta_0 k^{-t}$ and $\mu_k \gets \mu_1 k^{-t}$ with $t \gets 0.7$ for all $k \in \N{}$.  Estimates of $\kappa_{\nabla f}$, $L_{\nabla f}$, $\{\kappa_{c_i}\}_{i\in[m]}$, $\{L_{c_i}\}_{i\in[m]}$, $\{\kappa_{\nabla c_i}\}_{i\in[m]}$, and $\{L_{\nabla c_i}\}_{i\in[m]}$ were determined by randomly generating $n$ points from a normal distribution centered at $x_1$, then computing constraint function and derivative values at these points to compute bound and Lipschitz-constant estimates.  These values were employed to set the step sizes $\{\alpha_k\}$ as stated in Parameter Rule~\ref{par.deterministic} with $t_\alpha \gets 0$.  (This ensures $t + t_\alpha \in [-1,0)$, which means that the condition on these exponents in Theorems~\ref{th.deterministic} was satisfied.)  Finally, \texttt{SLIP} employed the parameter values $\eta \gets (\theta_0/\mu_1 + 1)/2$, $\underline\eta \gets \theta_0 + 10^{-8}$, $\underline\zeta \gets 1$, $\overline\zeta \gets 1$, and $\zeta \gets 1$.  Rather than attempt to ensure that \eqref{eq.d_conditions} holds for all $k \in \N{}$, \texttt{SLIP} follows the strategy described at the end of \S\ref{sec.assumption} to reset $\mu_1$ whenever a search direction is computed that does not satisfy~\eqref{eq.d_conditions}.  In particular, $\mu_1$ was multiplied by a factor of~2 in such cases, although we imposed an overall limit of $10^4$ beyond which $\mu_1$ was not increased.  One additional feature that we included in the implementation is that we allowed the algorithm to explore values of $\gamma_k$ greater than 1, since we believe this would be a useful feature in any practical implementation of the method.  Specifically, if $x_k + \alpha_k d_k \not\in \Ncal(\theta_k)$, then the algorithm would iteratively \emph{reduce} $\gamma_k$ from~1 until $x_k + \gamma_k \alpha_k d_k \in \Ncal(\theta_k)$ was found.  On the other hand, if $x_k + \alpha_k d_k \in \Ncal(\theta_k)$, then the algorithm would iteratively \emph{increase}~$\gamma_k$ from~1 as long as $x_k + \gamma_k \alpha_k d_k \in \Ncal(\theta_k)$ and the barrier term in the barrier-augmented objective function did not increase.

We ran \texttt{SLIP} with an iteration budget of $K := 2 \times 10^4$ in order to determine the progress that the algorithm could make within a budget.  (To resemble the stochastic setting, our implementation of \texttt{SLIP} does not employ second-order derivatives of the objective function, so one should not expect the algorithm to converge at a fast rate like other interior-point methods that employ second-order derivatives.)  Let $x_K \in \R{n}$ be the final iterate generated by a run of the algorithm on a test problem.  Since our algorithm is a feasible method, one can compare $f(x_1)$ and $f(x_K)$ for each run in order to determine if the algorithm has made progress.  We confirmed that, indeed, in all of our runs of \texttt{SLIP} over our entire test set, we found that $f(x_K) < f(x_1)$.  However, without knowledge of $f_{\inf} := \inf_{x \in \Fcal} f(x)$ for each problem in our test set, it is not possible to compare $f(x_K)$ with $f_{\inf}$.  Instead, we considered the relative stationarity measure (recall Theorem~\ref{th.deterministic}):
\bequation\label{eq.relstat}
  \text{relative stationarity} := \frac{\|P\nabla_x \phi(x_K,\mu_K)\|_2}{\min\{\|P\nabla_x \phi(x_1,\mu_1)\|_2, \|P\nabla_x \phi(x_1,\mu_K)\|_2\}}.
\eequation
Here, the $\min$ in the denominator respects the fact that different choices of the barrier parameter correspond to different multiplier estimates at the initial point (see \eqref{eq.multipliers}), which in turn give different stationarity measures at the initial point.  Our use of the $\min$ in the denominator in this definition of relative stationarity means that we are determining the better of the stationarity measures at the initial point to compare against our stationarity measure at the final iterate.

Figure~\ref{fig.cutest_histogram} shows a histogram of relative stationarity values over our test set of problems from the CUTEst collection.  The results show that, generally speaking, the relative stationarity values were quite small, meaning that the final iterates in most runs appeared to be much closer to stationarity than the initial iterates, thus indicating that within the iteration budget the algorithm made substantial progress toward stationarity.  We conjecture that the problems on which the relative stationarity value was higher may be ones that are very nonlinear and/or have constraints that are degenerate near the iterates generated by the algorithm.

\bfigure[ht]
  \centering
  \includegraphics[width=0.9\textwidth]{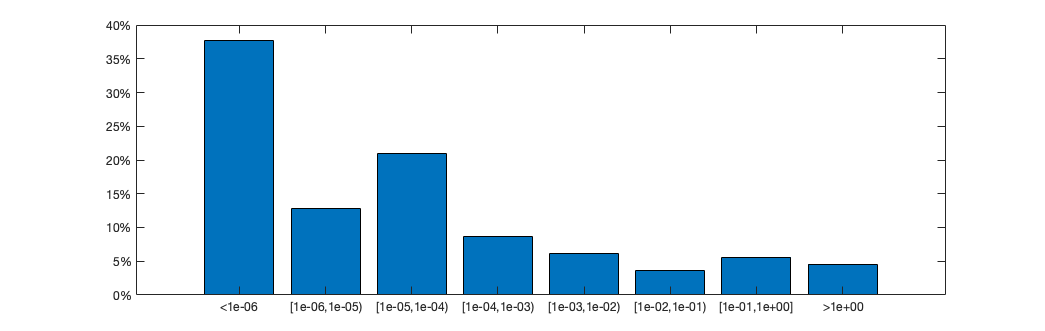}
  \caption{Histogram of relative stationarity values for experiments with a deterministic version of \texttt{SLIP} over problems from the CUTEst collection with a strictly feasible initial point.}
  \label{fig.cutest_histogram}
\efigure

As a second experiment, we tested stochastic runs of \texttt{SLIP} on a randomly generated instance of a second-order conic optimization problem, namely,
\bequation\label{prob.socp}
  \min_{x \in \R{n}}\ c^Tx\ \ \st\ \ Ax = b\ \ \text{and}\ \ \|x_{1:n-1}\|_2 \leq x_n.
\eequation
Specifically, we randomly generated $A$ and $b$ such that a strictly feasible point existed, and randomly generated $c$ such that the problem had a finite optimal value.  We chose a problem of this type since it is nonlinear, but convex, so for our first demonstration of stochastic runs of our algorithm we could avoid challenging situations that arise due to nonconvexity, such as the algorithm potentially converging to different minimizers.  (By contrast, our third experiment considers a nonconvex problem.)  First, as a benchmark, we ran deterministic \texttt{SLIP} with the same parameter choices as our experiments with the CUTEst collection, except for two differences: (a) we set $t_\alpha = -0.151$ so that $t + t_\alpha \in [-1,0)$ and $t + 2t_\alpha < -1$, as required in Parameter Rule~\ref{par.stochastic}, and (b) our procedure for computing $\gamma_k$ for all $k \in \N{}$ \emph{did~not} observe barrier-augmented objective function values, since this would require objective function values that are not tractable to compute in a stochastic setting.  Instead, we allowed $\gamma_k$ to increase above 1, but only up to a limit of 10, when doing so maintained all constraint values less than or equal to $-\theta_k$.

Second, we ran stochastic \texttt{SLIP} with 10 different seeds for the random number generator.  These runs used the same values for $x_1$, $\theta_0$, initial $\mu_1$, $\eta$, and $\underline\eta$ as the deterministic algorithm along with $\gamma_{k,\max} = 10$ for all $k \in \N{}$.  It also employed the same Lipschitz constant estimates and employed the same procedured for ``resetting'' $\mu_1$ (up to a limit of $10^4$) as the deterministic algorithm.  Thus, the main difference between the deterministic and stochastic runs of the algorithm were that the latter employed stochastic objective gradient estimates of the form $c + \Delta C$, where for all $k \in \N{}$ in each run the vector $\Delta C$ was drawn randomly from a standard normal distribution.  By design of the algorithm, the final iterates generated by the deterministic and all stochastic runs of the algorithm were feasible with respect to all of the constraints.  Table~\ref{tab.socp} shows the final objective values for each of the 10 stochastic runs of \texttt{SLIP}.  For reference, $f(x_1) =$ 1.80561e+03 for all runs, deterministic \texttt{SLIP} achieved $f(x_K) =$ 6.31100e+02, and the relative stationarity measure achieved by deterministic \texttt{SLIP} (see \eqref{eq.relstat}) was 1.65377e-02.  Thus, overall, one can observe that the deterministic and stochastic runs of \texttt{SLIP} yielded final iterates that were relatively very close to stationarity.

\begin{table}[ht]
  \centering
  \caption{Final objective function values over 10 runs of stochastic \texttt{SLIP} to solve \eqref{prob.socp}.}
  \label{tab.socp}
  \btabular{|c|ccccc|}
    \hline
    run & 1 & 2 & 3 & 4 & 5 \\
    \hline
    $f(x_K)$ & 6.31103e+02 & 6.31102e+02 & 6.31100e+02 & 6.31101e+02 & 6.31100e+02 \\
    \hline
    \hline
    run & 6 & 7 & 8 & 9 & 10 \\
    \hline
    $f(x_K)$ & 6.31101e+02 & 6.31102e+02 & 6.31101e+02 & 6.31101e+02 & 6.31103e+02 \\
    \hline
  \etabular
\end{table}

As a third experiment, we tested stochastic runs of \texttt{SLIP} to train a binary classifier subject to a norm constraint for the \texttt{mushrooms} data set from the LIBSVM collection \cite{CC01a}.  The model that we trained was an artificial neural network with a single hidden layer with 512 nodes, tanh activation, and cross-entropy loss.  This objective function is nonconvex.  We trained the model subject to the constraint that the neural network weights had squared $\ell_2$-norm less than or equal to 100.  This feasible region is convex.  Like for our second experiment, we chose a convex feasible region to ensure that it would be reasonable to expect multiple stochastic runs of the algorithm to converge to a point of similar solution quality in terms of the final objective function values (our measure of comparison), despite the nonconvexity of the objective function.  The experimental setup was the same as for our second experiment, except that, rather than artificial noise in the stochastic-gradient estimates, we employed mini-batch gradient estimates with a mini-batch size of 256.  (The \texttt{mushrooms} dataset has 8124 data points, each with 112 features.)

Table~\ref{tab.nn} shows the final objective values for each of the 10 stochastic runs of \texttt{SLIP}.  For reference, $f(x_1) =$ 6.93137e-01 for all runs, deterministic \texttt{SLIP} achieved $f(x_K) =$ 7.77697e-03, and the relative stationarity measure achieved by deterministic \texttt{SLIP} (see \eqref{eq.relstat}) was 7.06457610e-02.  Thus, one can observe that the deterministic and stochastic runs of \texttt{SLIP} yielded final iterates of comparable quality.

\begin{table}[ht]
  \centering
  \label{tab.nn}
  \btabular{|c|ccccc|}
    \hline
    run & 1 & 2 & 3 & 4 & 5 \\
    \hline
    $f(x_K)$ & 7.77711e-03 & 7.76221e-03 & 7.79686e-03 & 7.82982e-03 & 7.78385e-03 \\
    \hline
    \hline
    run & 6 & 7 & 8 & 9 & 10 \\
    \hline
    $f(x_K)$ & 7.79379e-03 & 7.80981e-03 & 7.80194e-03 & 7.75375e-03 & 7.77848e-03 \\
    \hline
  \etabular
\end{table}

%*********
% Section
%*********
\section{Conclusion}\label{sec.conclusion}

We have proposed, analyzed, and tested a single-loop interior-point framework for solving constrained optimization problems.  Of particular interest is a stochastic-gradient-based version of the framework, which can be employed for solving problems with affine equality constraints and (potentially nonconvex) nonlinear inequality constraints, at least when a strictly feasible initial point can be provided (or at least computed through a so-called Phase~I algorithm, such as the method that we present in Appendix~\ref{sec.appendix}).  We have shown that the algorithm possesses convergence guarantees in both the deterministic and stochastic settings.  We have also shown that a deterministic version of the algorithm performs reliably on a large test set of problems, and that a stochastic version of the algorithm yields good results both in the setting of artificial noise and in the context of a mini-batch stochastic-gradient-based algorithm for the training an artificial neural network for binary classification subject to a norm constraint.

A few significant open questions remain in the study of stochastic-gradient-based interior-point methods for solving nonlinearly constrained optimization problems.  For example, our theoretical convergence guarantees requires that the ratio $\mu_1/\theta_0$ is sufficiently large and that each computed search direction satisfies the conditions in Assumptions~\ref{ass.c} and/or \ref{ass.H}, but, as we have explained throughout the paper, a computationally effective strategy for computing search directions that adhere to our theoretical guarantees is not straightforward to design.  A related open question is whether it is possible to design a single-loop \emph{infeasible} interior-point framework that possesses theoretical convergence guarantees that are at least on par with those that we offer for the algorithm proposed in this paper.  The main challenge in the design of such an approach is that our theoretical guarantees, which employ a prescribed sequence $\{\mu_k\} \searrow 0$, rely heavily on the fact that the iterates generated by our algorithm are feasible in every iteration.  There exist stochastic-gradient-based infeasible Newton-type methods for solving equality-constrained optimization problems; see, e.g., \cite{BeraCurtRobiZhou21,BeraBollZhou22,BeraCurtOneiRobi21,BeraShiYiZhou22,CurtRobiZhou21,CurtOneiRobi21,FangNaMahoKola22,NaAnitKola22,NaMaho22,qiu2023sequential}.  Hence, one might expect it to be possible to employ such methods for solving the equality-constrained subproblems that arise in an infeasible interior-point method.  However, even if one were to introduce slack variables for which the barrier functions are introduced, it remains unclear how to merge such an approach with our strategy for ensuring feasibility at all iterates.  Concretely, suppose that inequality constraints $c(x) \leq 0$ are reformulated with slack variables as $c(x) + s = 0$ along with $s \geq 0$, where the latter inequalities are handled with a barrier function.  Attempting to follow the strategy for bound-constrained optimization in \cite{CKRW23}, one is confronted with the challenge that the search direction in the slack variables, call it $d_k^s$, depends on the search direction in the original variables, say $d_k^x$, where it may be possible that $c(x_k) \not\leq 0$.  It is difficult to follow a strategy such as ours to ensure, e.g., that there exists $\gamma$ \emph{uniformly bounded below} such that $s_k + \gamma \alpha_k d_k^s$ is sufficiently positive, say to stay within a prescribed neighborhood.  The difficulties of designing such an approach led us to propose the feasible method that has been presented in this paper, although it is possible that such an approach, or one based on alternative strategies, could be designed with convergence guarantees.

%**********
% Appendix
%**********
\appendix

\section{Appendix: Algorithm for Finding a Strictly Feasible Point}\label{sec.appendix}

Our algorithm for finding a strictly feasible initial point for our numerical experiments with problems from the CUTEst collection \cite{GoulOrbaToin15} is presented as Algorithm~\ref{alg.phaseI} in this appendix.  For simplicity, for our discussion here we state user-defined parameters in terms of the specific values that we used in our implementation rather than introduce a generic parameter range.

The algorithm can be understood as an infeasible interior-point method, where line searches on a merit function are performed as a step-acceptance mechanism.  Rather than attempt to solve an optimization problem to high accuracy, the aim of the algorithm is merely to find a strictly feasible point, i.e., a point in $\Fcal_{<0}$ (recall \eqref{eq.F_strict}) satisfying affine equality constraints and strictly satisfying (potentially nonlinear) inequality constraints.  Hence, the barrier parameter is fixed at 1 and the algorithm terminates as soon as a strictly feasible point is found.  A strictly feasible point is defined as follows.  As opposed to problem~\eqref{prob.opt}, which is stated in terms of only one-sided inequality constraints, our definition here of a strictly feasible point distinguishes between one- and two-sided inequalities.  With respect to any one-sided inequality, say, $\varphi(x) \leq \varphi_u$, strict feasibility of a point $x$ is defined as $\varphi(x) \leq \varphi_u - 10^{-4}$.  With respect to any two-sided inequality, say, $\varphi_l \leq \varphi(x) \leq \varphi_u$, strict feasibility of a point $x$ is defined as
\bequationNN
  \varphi_l + 10^{-4} \min\{\varphi_u - \varphi_l, 1\} \leq \varphi(x) \leq \varphi_u - 10^{-4} \min\{\varphi_u - \varphi_l, 1\}
\eequationNN

Supposing again that the inequality constraints are stated as in problem~\eqref{prob.opt}, the algorithm can be understood as an iterative method toward solving
\bequation\label{eq.phase_one_problem}
  \min_{(x,s) \in \R{n} \times \R{m}}\ -\sum_{i \in [m]} \log (s_i)\ \ \st \left\{ \baligned Ax &= b \\ c(x) + s &= 0 \ealigned \right.
\eequation
the first-order optimality conditions for which are
\bequationNN
  -S^{-1}e + z = 0,\ \ A^Ty + \nabla c(x)z = 0,\ \ c(x) + s = 0,\ \ \text{and}\ \ Ax = b.
\eequationNN
Given an initial point $x_0 \in \R{n}$, the algorithm commences by attempting to compute $x_1 \in \R{n}$ with $Ax_1 = b$ by solving a least-squares problem.  (If this yields $\|Ax_1 - b\|_2 \leq 10^{-6}$, then the algorithm continues; otherwise, the run is a failure.)  Then, up to an iteration limit, the algorithm follows a standard infeasible interior-point approach with line searches on a merit function; in particular, the merit function $\phi : \R{n} \times \R{m} \times \R{}_{>0} \to \R{}$ is defined by
\bequationNN
  \phi(x,s,\tau) = - \tau \sum_{i=1}^m \log (s_i) + \|c(x) + s\|_1,
\eequationNN
where $\tau \in \R{}_{>0}$ is a merit parameter that is updated adaptively.  (If the iteration limit is exceeded, then the run is a failure.)  Firstly, a Hessian approximation is determined with a modification, if necessary, to ensure that it is sufficiently positive definite.  To represent potential use in practice with our proposed Algorithm~\ref{alg.slip}, the Hessian approximation does not employ second-order derivatives of the objective function, although it does employ second-order derivatives of the constraint functions.  Secondly, a search direction is determined by solving a linear system.  Thirdly, the merit parameter is updated using a standard technique from the literature; see, e.g., \cite{ByrdHribNoce99}.  Specifically, defining the predicted reduction in the merit function corresponding to the tuple $(x_k,\tau_k)$ and search direction components $(d_k^x,d_k^s)$ as
\bequationNN
  \Delta q(x_k, s_k, d_k^x, d_k^s, \tau_k) = - \tau_k ( -{d_k^s}^T S_k^{-1} \ones + \thalf {d_k^x}^T H_k d_k^x + \thalf {d_k^s}^T S_k^{-2} d_k^s) + \|c(x_k) + s_k\|_1,
\eequationNN
the update ensures that $\tau_k$ is sufficiently small such that $\Delta q(x_k, s_k, d_k^x, d_k^s, \tau_k)$ is sufficiently large relative to $\|c(x_k) + s_k\|_1$.  Fourthly, the largest step size in $(0,1]$ satisfying a fraction-to-the-boundary rule is determined.  Finally, a line search is conducted to compute a step size that satisfies the fraction-to-the-boundary rule and yields sufficient decrease in the merit function.

\begin{algorithm}[ht]
  \caption{Finding a Strictly Feasible Point}
  \label{alg.phaseI}
  \begin{algorithmic}[1]
    \Require $x_0 \in \R{n}$, $\tau_0 = 1$
    \State Set $x_1 \gets \text{argmin}_{x \in \R{n}} \thalf \|x - x_0\|_2^2\ \st\ Ax = b$
    \State Set $s_1 \gets \max\{-c(x_1),1\}$ and $z_1 \gets s_1^{-1}$ (component-wise)
    \For{$k \in [10^3]$}
      \If{$c(x_k)$ is sufficiently interior, component-wise}
        \State terminate and return $x_k$
      \EndIf
      \State Set $H_k \gets I + \sum_{i\in[m]} [y_k]_i \nabla^2 c_i(x_k) + \lambda_k I$ for $\lambda_k \geq 0$ such that $H_k \succeq 10^{-4} I$
      \State Solve
      \begin{align*}
        \bbmatrix H_k & 0 & A^T & \nabla c(x_k) \\ 0 & S_k^{-2} & 0 & I \\ A & 0 & 0 & 0 \\ \nabla c(x_k)^T & I & 0 & 0 \ebmatrix \bbmatrix d_k^x \\ d_k^s \\ y_k \\ d_k^z \ebmatrix
        =
        - \bbmatrix \nabla c(x_k) z_k \\ -S_k^{-1}\ones + z_k \\ 0 \\ c(x_k) + s_k \ebmatrix 
      \end{align*}
      \State Set, with $\Theta_k := -{d_k^s}^T S_k^{-1} \ones + {d_k^x}^T H_k d_k^x + {d_k^s}^T S_k^{-2} d_k^s$, first
      \begin{align*}
        \tau_k^{\text{trial}} &\gets \begin{cases} \infty & \text{if $\Theta_k \leq 0$} \\ \tfrac{0.5 \|c(x_k) + s_k\|_1}{\Theta_k} & \text{otherwise} \end{cases} \\  
        \text{then}\ \tau_k &\gets \begin{cases} \tau_{k-1} & \text{if $\tau_{k-1} \leq \tau_k^{\text{trial}}$} \\ (1-10^{-6}) \tau_k^{\text{trial}} & \text{otherwise} \end{cases}
      \end{align*}
      \State Set $\alpha_k^{\text{ftb}} \gets \min\{ \alpha \in (0,1] : s_k + \alpha d_k^s \geq 0.1 s_k\}$
      \State Set $\alpha_k \gets \alpha^j \alpha_k^{\text{ftb}}$, where $j$ is the minimum value in $\{0\} \cup \N{}$ such that
      \bequationNN
        \phi(x_k + \alpha^j \alpha_k^{\text{ftb}} d_k^x, s_k + \alpha^j \alpha_k^{\text{ftb}} d_k^s, \tau_k) \leq \phi(x_k, s_k, \tau_k) - 10^{-4} \alpha^j \alpha_k^{\text{ftb}} \Delta q(x_k, s_k, d_k^x, d_k^s, \tau_k)
      \eequationNN
    \EndFor
  \end{algorithmic}
\end{algorithm}

%**************
% Bibliography
%**************
\bibliographystyle{plain}
\bibliography{ref}

\begin{thebibliography}{10}

\bibitem{BeraBollZhou22}
Albert~S. Berahas, Raghu Bollapragada, and Baoyu Zhou.
\newblock An adaptive sampling sequential quadratic programming method for
  equality constrained stochastic optimization.
\newblock {\em arXiv 2206.00712}, 2022.

\bibitem{BeraCurtOneiRobi21}
Albert~S. Berahas, Frank~E. Curtis, Michael~J. O'Neill, and Daniel~P. Robinson.
\newblock {A Stochastic Sequential Quadratic Optimization Algorithm for
  Nonlinear Equality Constrained Optimization with Rank-Deficient Jacobians}.
\newblock {\em Mathematics of Operations Research},
  https://doi.org/10.1287/moor.2021.0154, 2023.

\bibitem{BeraCurtRobiZhou21}
Albert~S. Berahas, Frank~E. Curtis, Daniel~P. Robinson, and Baoyu Zhou.
\newblock Sequential quadratic optimization for nonlinear equality constrained
  stochastic optimization.
\newblock {\em SIAM Journal on Optimization}, 31(2):1352--1379, 2021.

\bibitem{BeraShiYiZhou22}
Albert~S. Berahas, Jiahao Shi, Zihong Yi, and Baoyu Zhou.
\newblock Accelerating stochastic sequential quadratic programming for equality
  constrained optimization using predictive variance reduction.
\newblock {\em Computational Optimization and Applications}, 86:79--116, 2022.

\bibitem{ByrdHribNoce99}
Richard~H. Byrd, Mary~E. Hribar, and Jorge Nocedal.
\newblock An interior point algorithm for large-scale nonlinear programming.
\newblock {\em SIAM Journal on Optimization}, 9(4):877--900, 1999.

\bibitem{CC01a}
Chih-Chung Chang and Chih-Jen Lin.
\newblock {LIBSVM}: A library for support vector machines.
\newblock {\em ACM Transactions on Intelligent Systems and Technology},
  2:27:1--27:27, 2011.
\newblock Software available at \url{http://www.csie.ntu.edu.tw/~cjlin/libsvm}.

\bibitem{CKRW23}
Frank~E. Curtis, Vyacheslav Kungurtsev, Daniel~P. Robinson, and Qi~Wang.
\newblock A {{Stochastic-Gradient-based Interior-Point Algorithm}} for
  {{Solving Smooth Bound-Constrained Optimization Problems}}.
\newblock {\em arXiv e-prints}, arXiv:2304.14907, April 2023.

\bibitem{CurtOneiRobi21}
Frank~E. Curtis, Michael~J. O'Neill, and Daniel~P. Robinson.
\newblock {Worst-Case Complexity of an SQP Method for Nonlinear Equality
  Constrained Stochastic Optimization}.
\newblock {\em {Mathematical Programming}},
  https://doi.org/10.1007/s10107-023-01981-1, 2023.

\bibitem{CurtRobiZhou21}
Frank~E. Curtis, Daniel~P. Robinson, and Baoyu Zhou.
\newblock Inexact sequential quadratic optimization for minimizing a stochastic
  objective function subject to deterministic nonlinear equality constraints.
\newblock {\em arXiv 2107.03512 (to appear in INFORMS Journal on
  Optimization)}, 2021.

\bibitem{FangNaMahoKola22}
Yuchen Fang, Sen Na, Michael~W. Mahoney, and Mladen Kolar.
\newblock Fully stochastic trust-region sequential quadratic programming for
  equality-constrained optimization problems.
\newblock {\em SIAM Journal on Optimization}, 34(2):2007--2037, 2024.

\bibitem{PyCUTEst2022}
Jaroslav Fowkes, Lindon Roberts, and \'Arp\'ad Brmen.
\newblock Pycutest: an open source python package of optimization test
  problems.
\newblock {\em Journal of Open Source Software}, 7(78):4377, 2022.

\bibitem{GoulOrbaToin15}
Nicholas I.~M. Gould, Dominique Orban, and Philippe~L. Toint.
\newblock {CUTEst: a constrained and unconstrained testing environment with
  safe threads for mathematical optimization}.
\newblock {\em Comp. Opt. Appl.}, 60(3):545--557, 2015.

\bibitem{MangFrom67}
Olvi~L. Mangasarian and Stan Fromovitz.
\newblock {The Fritz John necessary optimality conditions in the presence of
  equality and inequality constraints}.
\newblock {\em Journal of Mathematical Analysis and Applications}, 17:37--47,
  1967.

\bibitem{NaAnitKola22}
Sen Na, Mihai Anitescu, and Mladen Kolar.
\newblock {An adaptive stochastic sequential quadratic programming with
  differentiable exact augmented Lagrangians}.
\newblock {\em Mathematical Programming}, 199:721--791, 2023.

\bibitem{NaMaho22}
Sen Na and Michael~W. Mahoney.
\newblock Asymptotic convergence rate and statistical inference for stochastic
  sequential quadratic programming.
\newblock {\em arXiv 2205.13687}, 2022.

\bibitem{NW99}
Jorge Nocedal and Stephen~J. Wright.
\newblock {\em Numerical Optimization}.
\newblock Springer Series in Operations Research and Financial Engineering.
  Springer-Verlag New York, 2006.

\bibitem{qiu2023sequential}
Songqiang Qiu and Vyacheslav Kungurtsev.
\newblock A sequential quadratic programming method for optimization with
  stochastic objective functions, deterministic inequality constraints and
  robust subproblems.
\newblock {\em arXiv 2302.07947}, 2023.

\bibitem{VandShan99}
Robert~J. Vanderbei and David~F. Shanno.
\newblock {An Interior-Point Algorithm for Nonconvex Nonlinear Programming}.
\newblock {\em Computational Optimization and Applications}, 13:231--252, 1999.

\bibitem{WaecBieg06}
Andreas W{\"a}chter and Lorenz~T. Biegler.
\newblock On the implementation of an interior-point filter line-search
  algorithm for large-scale nonlinear programming.
\newblock {\em Mathematical Programming}, 106(1):25--57, 2006.

\end{thebibliography}

\end{document}